\documentclass{amsart} \usepackage{amssymb,amsfonts,amscd}
\newtheorem{theorem}{Theorem}[section]
\newtheorem{lemma}[theorem]{Lemma}
\newtheorem{corollary}[theorem]{Corollary}
\newtheorem{proposition}[theorem]{Proposition}
\theoremstyle{remark}

\theoremstyle{definition}
\newtheorem{definition}[theorem]{Definition}

\numberwithin{equation}{section} \makeatother

\newtheorem*{xrem}{Remark}

\DeclareMathOperator{\Cdb}{{\mathbb C}}
\DeclareMathOperator{\Rdb}{{\mathbb R}}

\DeclareMathOperator{\Ndb}{{\mathbb N}}

\begin{document}

\title[Open projections in operator algebras I]{Open projections in operator algebras I: Comparison theory} \author[D. P. Blecher]{David P. Blecher}
\author[Matthew Neal]{Matthew Neal}

\address{Department of Mathematics, University of Houston, Houston, TX
77204-3008}
 \email[David P.
Blecher]{dblecher@math.uh.edu}
\address{Department of Mathematics,
Denison University, Granville, OH 43023}
\email{nealm@denison.edu}

\thanks{February 7, 2012}

\begin{abstract} We begin a program of generalizing
basic elements of the theory of
comparison, equivalence, and subequivalence, of elements in $C^*$-algebras, to the setting of
more general algebras.  In particular, we follow the recent lead of Lin, Ortega, R{\o}rdam, and Thiel
of studying these equivalences, etc.,   in terms of open projections or module isomorphisms.
We also define and characterize a new class of inner ideals in operator algebras, and  develop a matching theory of open partial isometries in operator ideals
which  simultaneously  generalize the open projections in operator algebras (in the sense of
 the authors and Hay),
and the open partial isometries (tripotents)
introduced by the authors.
\end{abstract}



\maketitle

\let\text=\mbox

\section{Introduction and notation}
Inspired by a recent paper of Ortega, R{\o}rdam, and Thiel \cite{ORT}, we begin a program of generalizing
basic elements of the theory of
comparison, equivalence, and subequivalence, of
elements in $C^*$-algebras, to the setting of
more general algebras.
To do this we establish several technical results and tools, such as
 a new class of inner ideals in operator algebras, and a
matching theory of `open partial isometries'
which  simultaneously generalize the open projections in operator algebras (in the sense of
 the authors and Hay \cite{BHN}),
and (essentially) the open partial isometries or tripotents introduced in  \cite{BN0}.

We begin by considering a relation between two elements $a$ and $b$
which one may define in any
monoid or algebra $A$: namely that there exists $x, y \in A$ with $a = xy, b = yx$.
If $A$ is a group, then this defines an equivalence relation.   In an algebra
this is not an equivalence relation in general.  In fact this fails to
be an equivalence relation even for the case $A = M_n$, the $n$ by $n$ matrices.
Indeed even in this simple case the
characterization of when two matrices $a$ and $b$ can be
related in this way is quite subtle \cite{JS}, and does not match any kind of
equivalence relation on a $C^*$-algebra that has been important in $C^*$-algebra
theory as far as we know.
   The question arises of how to fix
this problem in an operator algebra, by which we will mean a (not necessarily selfadjoint) subalgebra
of $B(H)$, for a Hilbert space $H$.  In a $C^*$-algebra $A$ the `fix' that works, on a large subset
of $A$, is
to insist that $y = x^*$ above; and then this defines an equivalence relation
$\sim$ on the positive cone $A_+$ of $A$.   This is sometimes
called {\em Pedersen equivalence}  (see e.g.\
 \cite{PedF}).
We will expand this `fix' to
a larger set than  $A_+$, and to more general operator algebras than $C^*$-algebras.

In $C^*$-algebra theory, coarser equivalence relations than Pedersen equivalence,
and matching notions of subequivalence or comparison, are becoming increasingly important \cite{Bla}.
For example, recently the study of Cuntz equivalence and subequivalence
has become one of the most important areas of $C^*$-algebra
theory (see e.g.\
\cite{APT,ORT}).    And {\em comparison theory} in
$C^*$-algebras, where one considers coarser notions
of ordering of elements in a $C^*$-algebra than the usual $\leq$, generalizing in some sense the crucially important comparison
theory of projections in a von Neumann algebra (which
is often most effectively done using traces),
is obviously important even in its own right.  Studying these equivalence relations and comparisons involves marshaling
a formidable array of tools, results, and perspectives.
Some of this is simply impossible for more general operator
algebras, due for example to scarcity of projections
  and hereditary subalgebras in general.   A {\em classification program} along
these lines would be misguided.
 Our goal in this project is simply to
transfer some small portion of these tools, results, and perspectives
to more general operator algebras than $C^*$-algebras.  Along the
way we were also led to develop some aspects of `noncommutative
topology and noncommutative function theory relative to an operator algebra'.
In the present paper we discuss
generalizations of Pedersen, Blackadar, and Peligrad-Zsid\'o equivalence, using the paper of Ortega, R{\o}rdam, and Thiel \cite{ORT} as our guide.
The main point is that the authors of \cite{ORT} recast the
 equivalences and subequivalences
mentioned above, and also Cuntz equivalence and subequivalence, in terms
of open projections (see also \cite{Lin}).
In view of our generalizations of
 open projections in earlier projects \cite{BHN,BN0,BRead},
it was tempting to follow their approach in a more general setting.

Our variants of Pedersen, Blackadar, and Peligrad-Zsid\'o equivalence lead us to introduce in Section 3
the notion of a $*$-{\em open tripotent}, which
may be viewed as a generalization of Akemann's open projections, and  a generalization of
our version of the latter for
general operator algebras from \cite{BHN}.   They also are essentially a `non-selfadjoint variant' of the
open tripotents which we introduced in \cite{BN0}.    These
 $*$-open tripotents play a central role in our equivalence and comparison theory.   One characterization of $*$-open tripotents is in
terms of {\em hereditary bimodules}, which encapsulate  the
algebraic structure of $\overline{aAb}$ when $a$ and $b$ are Pedersen or Blackadar equivalent.
We develop the theory of  hereditary bimodules in Section 4.
For example we show using results
from \cite{BRead} that every separable  hereditary bimodule is of the form $\overline{aAb}$
where $a$ and
$b$ are Pedersen equivalent; and that a general hereditary bimodule is the closure of an increasing union of
subspaces of the same form.    These sections are mostly of technical interest,
for later use in the development of the theory, however applications are given
there in the form of  rephrasings of our variant of Pedersen equivalence
in the language of $*$-open tripotents or hereditary bimodules.  In Section 5
we study a variant of `Blackadar equivalence and subequivalence'
of elements in our algebra
$A$, and prove the analogue of some results of Lin and Ortega-R{\o}rdam-Thiel \cite{Lin,ORT}.  For
example, elements $a$ and $b$ of the type we consider are Blackadar equivalent in our
sense iff $\overline{aA} \cong
\overline{bA}$ completely isometrically and as right $A$-modules.  It is also equivalent to
$\overline{aAb}$ being a ``principal'  hereditary bimodule, or to the support projections
of $a$ and $b$ being equivalent in an appropriate variant of the sense of Peligrad and Zsid\'o.
 In an appendix we
prove  TRO versions of a few of our earlier results.

We remark that much of our theory depends on the existence of $n$th roots of elements satisfying $\Vert 1 - a \Vert \leq 1$.  Thus we would expect that a certain
portion of this theory generalizes to bigger classes of
elements which have roots, and a smaller portion
may generalize further still to a bigger class of Banach
algebras using the facts about roots in Banach algebras
(see e.g.\ \cite{LRS}).

Blanket convention:
Throughout this paper, $A$ is a fixed operator algebra, and $B$ is
a $C^*$-algebra that contains it.  Sometimes $B$ will be generated by $A$;
that is there is no proper $C^*$-subalgebra of $B$ containing $A$.
The diagonal $\Delta(A) = \{ a \in A : a^* \in A \}$ is
a $C^*$-algebra and a subalgebra of $A$.
Concerning notation and background, which we will discuss next,
it will be helpful
for the reader to have easy access to several of the references, particularly \cite{BLM,BHN,BN0,BRead,BMP},
for example for more detail beyond what is presented here.

For us a {\em projection}
is always an orthogonal projection.
We recall that by a theorem due to  Ralf Meyer,
every operator algebra $A$ has a unique unitization $A^1$ (see
e.g.\ \cite[Section 2.1]{BLM}). Below $1$ always refers to
the identity of $A^1$ if $A$ has no identity.   We are mostly
interested in operator algebras with  contractive approximate
identities (cai's).  We also call these {\em approximately unital}
operator algebras.
If $A$ is a nonunital
operator algebra represented (completely) isometrically on a Hilbert
space $H$ then one may identify the unitization $A^1$ with $A + \Cdb I_H$.   The
second dual $A^{**}$ is also an operator algebra with its (unique)
Arens product, this is also the product inherited from the von Neumann
algebra $B^{**}$ if
$A$ is a subalgebra of a $C^*$-algebra $B$.  Meets and joins in
$B^{**}$ of projections in $A^{**}$ remain in $A^{**}$,
since these meets and joins may be computed in the biggest
von Neumann algebra contained inside $A^{**}$. Note that
$A$ has a cai iff $A^{**}$ has an identity $1_{A^{**}}$ of norm $1$,
and then $A^1$ is sometimes identified with $A + \Cdb 1_{A^{**}}$.

For sets $X$ and $Y$ we write $XY$ for the closure of the sum of products
of the form $xy$ for $x \in X, y \in Y$, and similarly for a product of three sets.
If $a \in A$ then we write $A_a$ for $\overline{aAa}$.
An {\em inner ideal} in $A$ is a subspace $D$ with $DAD \subset D$.
A {\em hereditary subalgebra} (HSA) of $A$ is an
inner ideal which has a cai.   For the theory of HSA's in general operator algebras
see \cite{BHN}.  These objects are in an order preserving,
bijective correspondence with the {\em open projections} $p \in A^{**}$, by which we mean that there
is a net $x_t \in A$ with $x_t = p x_t p \to p$ weak*.  These are
also the open projections $p$ in the sense of Akemann
\cite{Ake,Ake2} in $B^{**}$, where $B$ is a $C^*$-algebra containing $A$, such that
$p \in A^{\perp \perp}$.  Indeed the weak* limit of a cai for a HSA is
an open projection, and is called the {\em support projection} of the HSA.
Conversely, if $p$ is an open projection in $A^{**}$, then
 $pA^{**}p \cap A$ is a  HSA in $A$.  We write
$pA^{**}p \cap A$ as $_pA_p$, or simply as $A_p$.  Similarly,
if $q$ is another projection, then
$_pA_q = \{ a \in A : a = paq \}$.

We recall that a {\em closed projection} is the `perp' of
an open projection.   Suprema (resp.\ infima) of open
(resp.\ closed) projections  in $A^{**}$, remain in $A^{**}$,
by the fact mentioned two paragraphs earlier about meets and joins,
together with the $C^*$-algebraic case of these facts \cite{Ake,Ake2}.

In \cite{PZ}, Peligrad and Zsid\'o introduce a notion of equivalence for open projections
$p$ and $q$ in the bidual of a $C^*$-algebra $B$: We say $p$ and $q$ are
{\em Peligrad-Zsid\'o equivalent}, and write $p \sim_{\rm PZ} q$ if there is a
partial isometry $v \in B^{**}$
such that $p = v^* v, q = v v^*, v B_p \subset B, v^* B_q \subset B$.
Peligrad and Zsid\'o prove in \cite[Lemma 1.3]{PZ} that in fact $q$ being open and $v^* B_q \subset B$
are implied by the other conditions above.

We will need some background from \cite{BN0}.  A {\em ternary ring of
operators} (or {\em TRO} for short), is a closed
subspace $Z$ of a C*-algebra $A$ such that $Z Z^* Z \subset Z$.
A {\em tripotent} is an element $u \in Z$ such that $u u^* u = u$.  This is clearly simply a partial isometry in $Z$.
We order tripotents by $u \leq v$ if and only if $u v^* u = u$.
This turns out to be equivalent to $u = v u^* u$, or to $u =  u u^* v$,
and implies that $u^* u \leq v^* v$ and  $u u^* \leq v v^*$.
The linking C*-algebra $L(Z)$ of a TRO $Z$  has `four corners'
$Z Z^*$, $Z, Z^*,$ and $Z^* Z$.   Here $Z Z^*$ is the closure of the
linear span of products $z w^*$ with $z, w \in Z$, and similarly
for $Z^* Z$.  The second dual of a TRO $Z$ is a TRO,
which is studied in terms of the von Neumann algebra
which is the second dual of $L(Z)$.
An {\em inner ideal} (resp.\ {\em ternary ideal}) of a TRO $Z$
is defined to be a closed subspace $J$ with $J Z^* J \subset J$
(resp.\ $J Z^* Z \subset J$ and $Z Z^* J \subset J$).

In most of our paper, although we will be using
notation and concepts from the theory of TRO's and JB*-triples,
in fact the TRO concerned will simply be the $C^*$-algebra $B$.
Thus we will write the discussion in the next few paragraphs
in terms of $B$, although it all makes sense even if $B$ were a
TRO.   The {\em Peirce $2$-space} of a tripotent $u$ in
 $B$ is the subset
$$B_2(u) = \{ z \in B :
z = u u^* z u^* u \} =  u u^* B u^* u.$$ Clearly $B_2(u)$ is an
inner ideal of $B$ in the TRO sense. There is a natural product (the
{\em Peirce product} $x \cdot y = x u^* y$) and involution (namely
$x^\sharp =  u x^* u$) on $B_2(u)$ making the latter space into a
unital $C^*$-algebra.  The identity element is $u$.  It is easy to
check that $u^* B_2(u)$ is a  C*-subalgebra of $B$ (this will be
$B^{\ast} B$ if $B$ is a TRO), and the map $z \mapsto u^* z$ is a
$*$-isomorphism from $B_2(u)$, with the product and involution
above, onto this C*-subalgebra.  We will write $B(u)$ for
$B^{**}_{2}(u) \cap B
 = \{ z \in B : z =
u u^* z u^* u \}$.  Here $B^{**}_{2}(u)$ is the
Peirce $2$-space of $u$ in $B^{**}$.

We will also occasionally use
other facts from Sections 2 and 3 of
\cite{BN0}, like those concerning the
{\em range tripotent} $r(x)$ of an element
 $x \in B$.  Since this is simply the
partial isometry occurring in the polar decomposition of $x$,
 so $x = r(x) |x|$, the facts concerning range tripotents are fairly simple
 and essentially well known.  For example, $r(x) \in B^{**}$ and $x = |x^*| r(x)$.

In \cite{BN0} we defined
 a tripotent $u$ in the second dual of a TRO
to be {\em open}, if when we consider
the Peirce $2$-space for $u$ as a W*-algebra in the  way described above,
then $u$ is the weak* limit in the second dual, of an increasing net of
positive elements from $B(u)$ (positive
with respect to the $C^*$-algebra
structure determined by the Pierce product;
and $B(u)$ is replaced by $Z(u)$ in the case of a general TRO
$Z$ rather than $B$).  Beware that this definition differs
from the one given in  \cite{ER3}. For example,
all unitaries are open in the sense of
that paper.   See \cite{FP,FP2,EF} for a recent JB*-triple
generalization of our notion.

There are several equivalent definitions
of open tripotents in
\cite[Theorem 2.10, Corollary 3.4]{BN0}.   For example,
a tripotent $u$ in the second dual is open iff
it is a weak* limit of an increasing net of range tripotents (defined
above), or if the projection
$$\hat{u} = \frac{1}{2} \left[ \begin{array}{ccl} u u^* & u \\
u^* & u^* u \end{array} \right]$$
is an open projection in Akemann's sense (discussed
earlier).   Another important characterization of
open tripotents, as the `support tripotent' of
certain inner ideals, will be
discussed in Section \ref{hbm}.  In  Proposition \ref{pzis}, we show that
our open tripotents are (at least in the case that the
TRO is a $C^*$-algebra) the partial isometries occurring in
the Peligrad and Zsid\'o equivalence of open projections defined above.   A special case of
this was proven as the equivalence of (i) and (vi) in \cite{BW}.
In our work \cite{BN0,BW} we were unaware of \cite{PZ}, which in retrospect clearly
 overlaps in small part with \cite{BN0}.  Namely, as we
have just seen,
 \cite{PZ} implicitly contains the notion of an open tripotent;
 and a couple of the
 results in \cite{PZ} may thus
be read as results about open
tripotents (these results were not repeated in
 \cite{BN0}).

In what follows we use the notation $\mathfrak{S}_A$ for either $\{
a \in A : \Vert 1 - 2a \Vert \leq 1 \}$, or $\{ a \in {\rm Ball}(A)
: \Vert 1 - a \Vert \leq 1 \}$.  Of course the first of these two
sets is contained in the second.  However the reader should choose
which of these two they prefer, and stick with this choice for the
rest of the paper. Note that $\mathfrak{S}_A$ is convex, and closed
under multiplication by scalars in $[0,1]$.   We write
$\mathfrak{c}_A$ for the cone $\Rdb^+ \mathfrak{S}_A$. This cone
will, just as in \cite{BRead},
 play a role for us
very much akin to the role of the positive cone in a $C^*$-algebra.  Indeed an underlying philosophy
in \cite{BRead} and the present paper and its sequel \cite{BN2} is to use this
cone to generalize important facts and theories for $C^*$-algebras which employ the positive cone, to more general operator algebras.  The set $\mathfrak{c}_A$ is `large': for example in a unital operator algebra
it includes $\Rdb^+(1 + {\rm Ball}(A))$, which spans $A$.  On the other hand, if $A$ has a cai then that cai may be chosen in this cone \cite{Read},
and $A = \overline{{\rm Span}}({\mathfrak S}_A)$.  To see
the latter, note that by \cite[Lemma 8.1]{BRead}
and the unital case just discussed, the weak* closure $E$ of $\overline{{\rm Span}}(\mathfrak{S}_A)$ equals $A^{**}$, hence by the bipolar theorem $A = \overline{{\rm Span}}(\mathfrak{S}_A)$.

We will use throughout the fact from  \cite{BRead} that elements in
$\mathfrak{S}_A$ have $n$th roots for all $n \in \Ndb$,
which are again in $\mathfrak{S}_A$.
There is an  explicit formula
for these roots, and indeed
for $x^t$ for all $0 < t \leq 1$
and $\Vert 1 - x \Vert \leq 1$ (the series in  Lemma \ref{unir} below).
If $a \in \mathfrak{S}_A$, then $(a^{\frac{1}{n}})$ converges
weak* to an open projection which is written as
$p_a$ or $s(a)$, and this is both the left
and the right support projection of $a$ (see \cite[Section 2]{BRead}).   If $a \in \mathfrak{c}_A$,
its support projection $s(a) = p_a$ is $s(\frac{a}{\Vert a \Vert})$.
In this case $\overline{aAa}$ is a HSA of $A$, and the
 support projection of this HSA is $s(a)$,
 so $\overline{aAa} = A_{p_a} = p_a A^{**} p_a
 \cap A$.

We will need some simple properties of the roots mentioned above:

\begin{lemma} \label{pow}  If $a \in  {\mathfrak S}_A$
then $a^{rs} = (a^r)^s$ if $r,s \in (0,1]$, and
$a^r a^s = a^{r+s}$ if $r,s, r+s \in (0,1]$.
\end{lemma}

\begin{proof}  This follows from the disk algebra
functional calculus used to define these powers in \cite{BRead}.
\end{proof}

\begin{lemma} \label{unir}  If $a \in  {\mathfrak S}_A$
then $a$ has a unique square root in ${\mathfrak S}_A$.
Indeed for $n \in \Ndb$, $a$ has a unique $n$th root
 in ${\mathfrak S}_A$ whose numerical range
is contained in the set of numbers $r e^{i \theta}$ with $|\theta| \leq \frac{\pi}{n}$.
This $n$th root is $\sum_{k = 0}^\infty \,
{t \choose k} (-1)^k (1-a)^k$, where $t = \frac{1}{n}$, and this
is a norm limit of polynomials in $a$ with no constant term.
\end{lemma}

\begin{proof}   This follows from a well known operator theoretic fact concerning roots (see e.g.\
\cite[Theorem 0.1]{LRS}), together with \cite[Proposition 2.3]{BRead}.
 We remark that
an early version of \cite{BRead} contained this
explicit formula, the published version referred to
$\sum_{k = 0}^\infty \,c_k  (1-x)^k$ without  explicitly
stating the values of $c_k$.  However that $c_k = {t \choose k} (-1)^k$
is obvious.   \end{proof}

\begin{xrem} A clarification, for our operators $a$,
one cannot expect $(a^{\frac{1}{2}})^* a^{\frac{1}{2}}
= |a|$.    In other words,
our `square root' $f(a) = a^{\frac{1}{2}}$ does
not satisfy $f(a^* a) = f(a^*) f(a)$.  To see
this note that if $a$ is any invertible matrix in
${\mathfrak S}_{M_n}$, then  $b =
a^{\frac{1}{2}}$ is invertible, and so
if $b^* b b^* b = a^* a = (b^*)^2 b^2$ then
 $b b^* = b^* b$.  So $b$ and hence $a$ is normal.
But this is false in general. \end{xrem}

\begin{lemma} \label{cofr}  If $S \in B(H)$ with
$\Vert I_H - S \Vert \leq 1$, and if $T$ is a
contraction in $B(K,H)$, then
$\Vert I_K - T^* S T \Vert \leq 1$.
\end{lemma}

\begin{proof} $I_K - T^* S T = V^* D V$, where
$V^* = [ (I_K - T^* T)^{\frac{1}{2}} \; \; \; T ]$
and $D = {\rm diag}(I_K , I_H - S)$.   Since $V$ is
an isometry the result is clear.
\end{proof}

We note that if $a  \in {\mathfrak c}_A$
and $n \in \Ndb$, then it is
clear that the left and right supports of $a^n$
equals $s(a)$.

\begin{lemma} \label{vpow}  If $a \in {\mathfrak S}_A$
and $v$ is a partial isometry in any
 containing $C^*$-algebra $B$ with $v^* v = s(a)$, then
$v a v^* \in {\mathfrak S}_B$ and
$(v a v^*)^r = v a^r v^*$ if $r \in (0,1) \cup \Ndb$.
\end{lemma}

\begin{proof}   This is clear if $r  = k \in \Ndb$.
Also $v (v^* v -a)^k v^* = (v v^* - v a v^*)^k$.
  For $r \in (0,1)$ we have $v a^r v^*$ equal to
 \begin{align*} \sum_{k = 0}^\infty \,
{r \choose k} (-1)^k v (1-a)^k v^* &= \sum_{k = 0}^\infty \, {r
\choose k} (-1)^k (v v^* - v a v^*)^k \\ &= \sum_{k = 0}^\infty \,
{r \choose k} (-1)^k (1- v a v^*)^k, \end{align*} which equals $(v a
v^*)^r,$ using Lemma \ref{cofr}.
\end{proof}

The following result will not be used here, but is of independent interest.  It may be generalized
to give the variant of \cite[Lemma 3.12]{Lin} appropriate to the setting of quotients of `rigged modules' in the sense of \cite{AGOH}.

\begin{proposition} \label{liq}  If $J$ is a right ideal in $A$ with a left cai $(e_t)$ in ${\mathfrak S}_A$,
for example if $J = \overline{aA}$ for some $a \in  {\mathfrak S}_A$, and if
$x \in A$ then the norm of $x+J$ in $A/J$ equals $\lim_t \Vert (1-e_t) x \Vert$.
 \end{proposition} \begin{proof}   Just as in the proof of  \cite[Lemma 3.12]{Lin}.  \end{proof}

We now consider $a, b \in {\mathfrak S}_A$, and are interested in
$\overline{aAb}$, an inner ideal in $A$.

\begin{proposition} \label{aab}  If $a, b \in {\mathfrak S}_A$, and
 $p = p_a, q = p_b$,
then $$\overline{aAb} = \, _pA_q = \{ x \in A : x = paq \} =
\overline{aA} \cap \overline{Ab} = \overline{aA} \, \overline{Ab} =
A_a A A_b.$$
 \end{proposition} \begin{proof}  Clearly $aAb \subset \, _pA_q$.
On the other hand if $x \in \,  _pA_q = \, _pA \cap  A_q =
\overline{aA} \cap \overline{Ab}$ then $a^{\frac{1}{n}} x
b^{\frac{1}{n}} \to x$.  So $x \in \overline{aAb}$. So
$$\overline{aAb} = \, _pA_q = \{ x \in A : x = paq \} =
\overline{aA} \cap \overline{Ab} \supset \overline{aA} \,
\overline{Ab} .$$ But of course $a A b = (a A b^{\frac{1}{2}})
b^{\frac{1}{2}} \subset \overline{aA} \, \overline{Ab},$ and so
$\overline{aAb} = \overline{aA} \, \overline{Ab}$. Clearly $A_a A
A_b \subset \overline{aAb}$, and so $A_a A A_b = \overline{aAb}$
since $a \in A_a, b \in A_b$.   \end{proof}

\section{A Pedersen type equivalence in operator algebras}

 In order to fix the problem with the relation
considered in the second paragraph of our paper, it is tempting
to define a relation $a \sim_c b$ if there exist
$x, y \in {\rm Ball}(B)$ with
 $a = x y, yx = b$.
As we will see momentarily this is still not what one would want, but
nonetheless this notation $\sim_c$ will be useful.  Indeed
$\sim_c$ is not
an equivalence relation on ${\mathfrak S}_A$. However it is almost an
equivalence relation in the following way:
$a \sim_c b$ and $b \sim_c d$ implies that $a^2 \sim_c d^2$.  Indeed if $a = x y, yx = b = wz, d = zw$ then $a^2 = xy xy = x b y = x wz y$
and $zy xw = z b w = z wz w = d^2$.

The relation $\sim_c$ is
certainly not a definition that we will consider
seriously, even if one insists on extra conditions on the
support projections for $x$ and $y$.  Indeed even in the simple case $A = M_2$
this definition
would be wrong in our context.  It differs completely from our simplest
equivalence relation, Pedersen equivalence, which in this case is just
unitary equivalence.  To see this,
let $x = {\rm diag} \{ \frac{1}{\sqrt{2}} ,
  \frac{1+ 3K}{\sqrt{2}} \}$, for some small positive $K$ to be determined,
and let $y$ be the $2$ by $2$ matrix with rows $\frac{1+K}{\sqrt{2}}, \sqrt{2} K$
and $0, \frac{1}{\sqrt{2}}$.  Then $xy$ is the $2$ by $2$ matrix with rows $\frac{1+K}{2},
K$ and $0, \frac{1+ 3K}{2}$, whereas $yx$ is the $2$ by $2$ matrix with rows $\frac{1+K}{2},
(1+3K)K$ and $0, \frac{1+ 3K}{2}$.   For $K$ small enough both $xy$ and $yx$
are in ${\mathfrak S}_A$ (indeed $x$ and $y$ are in   ${\mathfrak S}_A$  too).
 Simple formulae for the norm of an
upper triangular $2$ by $2$ matrix show that $xy$ and $yx$ have different
norms in general, hence are not unitarily equivalent.

The following seems to be a much more promising relation:

\begin{definition}  We say that $a$ and $b$ are {\em root equivalent}, and write $a \sim_r b$, if $a^{\frac{1}{n}} \sim_c b^{\frac{1}{n}}$
for all $n \in \Ndb$.  This is not
the same as $a \sim_c b$ on ${\mathfrak S}_A$, even in the case $A = M_2$.
This may be seen to be an equivalence relation on ${\mathfrak S}_A$ using the
fact at the end of the first paragraph of this section.  \end{definition}

  \begin{xrem}  (1) \ If $a, b \in  {\mathfrak S}_A$,
with $a \sim_r b$,  consider the four-tuple or `context'
$(aAa,bAb,aAb,bAa)$. If $a^{\frac{1}{3n}} = x y , b^{\frac{1}{3n}}
=yx$ for some $x, y \in {\rm Ball}(A)$, then letting $x' = xyx, y'=
yxy$, we have $a^{\frac{1}{n}} = x' y'$ and $b^{\frac{1}{n}} = y'
x'$, and $x' \in \overline{aAb}$ and $y' \in \overline{bAa}$. Since
these roots of $a$ and $b$ constitute cai's
 for $A_{a}$ and $A_{b}$ respectively,  $(\overline{aAa},
\overline{bAb}, \overline{aAb}, \overline{bAa})$ is a Morita context in the sense of
\cite{BMP}.

(2) \   Suppose that $a^{\frac{1}{n}} \sim_c b^{\frac{1}{n}}$ for {\em some} $n \in \Ndb$,
so $a^{\frac{1}{n}}  = x y, yx = b^{\frac{1}{n}}$.
Then $a^{\frac{1}{n}} x
= x y x = x b^{\frac{1}{n}}$.  Thus $a x =
(a^{\frac{1}{n}})^n x
= x (b^{\frac{1}{n}})^n = x b$, and similarly $a^m x = x b^m$
for all $m \in \Ndb$.  So $f(a) x = x f(b)$ for any function $f$ in the disk algebra, and even for
 more general functions.  In particular, $p_a x = x p_b$.  Similar assertions apply to $y$.
\end{xrem}

Parts of the following theorem use operations in
the containing $C^*$-algebra $B$.  However
it will be seen that it does not matter
which containing $C^*$-algebra we use, since other
parts such as (iii) do not use $B$ at all.
The careful reader will also notice apparent redundancy in some of the
following statements and proofs; however this is necessary
to obtain the first part
of the last assertion of the theorem (about being equivalent
with the same $x$ and $y$).

\begin{theorem}  \label{varp}  Suppose that
$a,b \in {\mathfrak S}_A$, and let $c = a^{\frac{1}{2}}$ and $d =
b^{\frac{1}{2}}$.  TFAE: \begin{enumerate}
\item [{\rm (i)}]  There exist
$w,y \in A$ with $a = wy, b = yw$ and $|y| = |c|$.
\item [{\rm (ii)}]  There exist
$x,y \in A$ with $a = xy, b = yx$ and $|y| = |c|$, and $x = x p_b$.
\item [{\rm (ii)}$^{\prime}$]  There exist
$x,y \in A$ with $a = xy, b = yx$ and $|y| = |c|$,  and $|x^*| =
|c^*|$.
\item [{\rm (ii)}$^{\prime \prime}$]  There exist
$x,y \in A$ with $a = xy, b = yx$ and $$|y| = |c|, |y^*| = |d^*|, |x|
= |d|, |x^*| = |c^*|.$$
\item [{\rm (ii)}$^{\prime \prime
 \prime}$]  There exist
$x,y \in A$ with $a = xy, b = yx$
and $x = c R$ and $y = S c$ for some contractions $R,S$
in $A^{\ast\ast}$ (or if one prefers, in $B^{\ast\ast}$).
 \item [{\rm (iii)}]  For all $n \in \Ndb$, there exist
$x_n,y_n \in {\rm Ball}(A)$ with $a^{\frac{1}{n}} = x_n y_n,
b^{\frac{1}{n}} = y_n x_n$, and  the sequence $( y_n a)$ has a norm
convergent subsequence.
    \item [{\rm (iv)}]  There exists a
    partial isometry $v \in B^{**}$ with
$p_a = v^* v, p_b = v v^*$, and  $v a,
   a v^* \in A$ and $b = v a v^*$.
\item [{\rm (iv)}$^{\prime}$]  There exists
$v \in \Delta(A^{**})$ with  $p_a = v^* v,$
and  $v a \in A$, and $b = v a v^*$.
\end{enumerate}
These imply that $\Vert a \Vert = \Vert b \Vert$.
The partial isometry $v$ in {\rm (iv)} above
may be chosen to be open and in $\Delta(A^{\perp \perp})$.
Also, the clauses {\rm (ii),(ii)$^{\prime}$,(ii)$^{\prime \prime}$,
(ii)$^{\prime \prime
 \prime}$} are equivalent to each other with the same $x$ and $y$, and these
imply that $x \in \overline{aAb}, y \in \overline{bAa}$.
\end{theorem}

\begin{proof} Note that (iv) clearly implies that
$\Vert a \Vert = \Vert b \Vert$.

(ii)$^{\prime \prime
 \prime}$  $\Rightarrow$ (ii)$^{\prime}$ \
Given (ii)$^{\prime \prime
 \prime}$, $a = xy = c R S c$,  we may assume that
$S = S  p_a$ and $p_a R = R$.   So
$cRS a = c a$ and so $c RS  = c RS  p_a = c$.
Thus $RS = p_a RS = p_a$.   It follows by e.g.\ \cite{Mb}
that $S$ is a
partial isometry and $R = S^*$, so that
$y^* y  = c^* S^* S c = c^* p_a c = c^* c$.  Similarly,
$x x^* = c c^*$.

(ii)$^{\prime}$ $\Rightarrow$ (ii)$^{\prime \prime
 \prime}$  \ Obvious.

(ii) $\Rightarrow$ (ii)$^{\prime \prime}$ \  Assume that there exist
$x,y \in A$ with $a = xy, b = yx$ and
$|y| = |c|$.
 Then $y = r(y) |c| =
v c$, where $v = r(y) r(c)^*$.  Hence $y = y p_a$, so that $r(y)^* r(y) \leq p_a$.  However
 $r(y)^* r(y) \geq p_a$ since $a = xy$, so $r(y)^* r(y) = p_a$.
Similarly, $r(y) r(y)^* \geq p_b$.  Note that
$v^* v = r(c) p_a r(c)^* = p_a$.  So $v$ is a tripotent.  Also,
$v v^* = r(y) p_a r(y)^* = r(y) r(y)^* \geq p_b$,
and we will see momentarily that this is an equality.

Note that $xy = x v c = c^2$ implies $xv a = ca$,
and hence $xv = xv p_a = c p_a = c$.
Since $y = vc$ we have   $va = v c^2 =
y c = y x v  = bv$.
Hence $y = vc = dv$.
It is now clear
that $p_b y = y$, so that $p_b \geq r(y) r(y)^*$,
giving $p_b  = r(y) r(y)^* =  v v^*$.
 Also $y y^* = d vv^* d^*
= d p_b d^* = d d^*$.

We now see that
$x = c v^*$, assuming that $x = x p_b$ (since then
$c v^* = x v v^* = x p_b = x$).
Thus $x x^* = c c^*$.  By symmetry to the case for $y$,
we must have $x^* x = c^* c$.

(ii) $\Rightarrow$ (iv) \ We use facts from the
proof above, noting that
$va = vc^2 = yc \in A, av^* = c^2 v^* = c x \in A$,
and $v a = v c^2 = d^2 v = bv$, so that $b = v a v^*$.

We now prove the assertions at the end of the theorem,
starting from (ii), and using the notation and facts
already established above.
To see that $v$ is open, consider $R_{r(c)^*}$, right
multiplication by  $r(c)^*$, on the TRO $\langle y \rangle$ generated by $y$.  This is a one-to-one ternary morphism,
since $y r(c)^* r(c) = yp_a = y$.  We show it also
maps into $B$:  note  $y r(c)^* c = y|c| \in B$,
so $y r(c)^* a \in B$, hence $y r(c)^* a^{\frac{1}{n}}
\in y r(c)^* \overline{aA} \subset B$.
  However $y r(c)^* = r(y) |c| r(c)^* = r(y)c^*$.
Now $c^* (a^{\frac{1}{n}})^* \to c^*$, so
$c^* a^{\frac{1}{n}}  \to c^*$ by 2.1.6 in \cite{BLM}.
Hence $y r(c)^* = r(y) c^* = \lim_n r(y)c^* a^{\frac{1}{n}}
\in B$.   So $R_{r(c)^*}(\langle y \rangle) \subset B$.    Hence $v$, the image under this ternary morphism of $r(y)$, is open, since
$r(y)$ is open.

As we saw two paragraphs back, $va \in A$ and $av^* \in A$.
By Lemma \ref{unir} we have
$v a^{\frac{1}{n}} \in A$.  In the weak* limit,
$v = v p_a \in A^{**}$.  Similarly, $a^{\frac{1}{n}} v^* \in A$ and so
 $v^* = p_a v^* \in A^{**}$.
So $v \in \Delta(A^{**})$.

By Proposition \ref{aab}, $x \in \, _{p_a}A_{p_b} = \overline{aAb}$.
Similarly, $y \in \overline{bAa}$.

(ii)$^{\prime}$ $\Rightarrow$ (ii) \ As in the proof  above,
we have $v = r(y) r(c)^*$ is a tripotent
with $y = vc, v^* v = p_a, v v^* = r(y) r(y)^* = p_b$.
  By symmetry, $w = r(c)^* r(x)$ is a tripotent
with $x = cw, w w^* = p_a, w^* w = r(x)^* r(x) = p_b$.
 Now $a = xy = cwvc$, and by the idea in the proof of
(ii)$^{\prime \prime \prime}$ $\Rightarrow$ (ii)$^{\prime}$ above,    it follows that $w = v^*$, so that
$x = c v^* = c v^* p_b = x p_b$.

(iv) $\Rightarrow$ (ii)$^{\prime}$ \ Set $x = c v^*, y = vc$.  These
are in $A$ by Lemma \ref{unir}.
  The rest
is easy from facts in the proof of the last two implications.

(i) $\Rightarrow$ (ii)$^{\prime}$ \ Let $x = wyw, y' = ywy$.
Then $|y'|^2 = a^* y^* y a = |a^{\frac{3}{2}}|^2$.
The proof that (ii) implies (ii)$^{\prime \prime}$ and (iv)  above,
with $a,b$ replaced by $a^3, b^3$,
 yields that $x = a^{\frac{3}{2}} v, y' = v a^{\frac{3}{2}}$
for a partial isometry $v$ satisfying
$v^{*}v=s(a)=s(a^{3})$, $v^{*}v=  s(b) = s(b^{3})$, and $va^{3}=b^{3}v$,
which implies that
$vc=dv$.
 Since $c \in \overline{aA}$ and
 $aA = ccA \subset c \overline{aA}
\subset  \overline{c^3 A}$,
$$vc \in  v \overline{aA} \subset v \overline{c^3 A}
\subset  \overline{v c^3 A} = \overline{y' A} \subset A .$$
 Similarly, $c v^* \in A$.  Then $a = c v^* vc,
b = d v v^* d =  vc^2 v^*$ as desired in  (ii)$^{\prime}$.

(ii) $\Rightarrow$ (i) \ Trivial.

(ii) $\Rightarrow$ (iii) \ Using the notation and proof above,
if $y_n = v a^{\frac{1}{2n}}, x_n =  a^{\frac{1}{2n}} v^*$,
then these are in $A$ by considerations used above,
and $x_n y_n = a^{\frac{1}{n}}$.  As we said above,
$v a = bv$, hence $v a^n = b^n v$, and so $v a^{\frac{1}{n}}
= b^{\frac{1}{n}} v$ by
considering the power series in Lemma \ref{unir}.
Hence $y_n x_n = b^{\frac{1}{2n}} v v^* b^{\frac{1}{2n}} = b^{\frac{1}{n}}$.

(iii) $\Rightarrow$ (iv)$^{\prime}$ \
Note that $y_n a^{\frac{1}{n}} = y_n x_n y_n = b^{\frac{1}{n}} y_n$
so that $y_n a = b y_n$. Similarly, $a x_n = x_n b$.
Replacing $y_n$ by $y_n x_n y_n$ and
$x_n$ by $x_n y_n x_n$ as before, but still calling
them $y_n$ and $x_n$, we have
$a^{\frac{3}{n}} = x_n y_n$ where
$x_n = p_a x_n, y_n = y_n p_a$, and
we still have $y_n a = b y_n$ and $a x_n = x_n b$, and
$(y_n a)$ still has a  convergent subsequence.
Suppose that $x_{n_k} \to S$ and $y_{n_k} \to T$ weak*, and $y_{n_k} a \to Ta$ in norm.
We have $S = p_a S$ and $T = T p_a$.
Then  $x_{n_k} y_{n_k} a
\to S T a = p_a a = a$.  As in the proof of
(ii)$^{\prime \prime \prime}$ $\Rightarrow$ (ii)$^{\prime}$,
this yields  $ST = ST p_a = p_a$, and $S^* = T$.
 Thus $T$ is a partial isometry $v \in \Delta(A^{**})$,
with $v^* v = p_a$.
Since $y_{n_k} a \to Ta$ in norm, $va \in A$.
Since  $y_n a = b y_n$ and $a x_n = x_n b$,
it follows that $b y_{n_k} x_{n_k} = y_{n_k} a x_{n_k}
\to v a v^*$, so that $b  = v av^*$.
Hence $p_b \leq v v^*$, and in fact
they are equal because $p_b v v^* = v v^*$ since
$p_b y_n = y_n$.
Thus we have (iv)$^{\prime}$.

(iv)$^{\prime}$ $\Rightarrow$ (iv) \ By Lemma \ref{vpow},
 $b^{\frac{1}{n}} = (v a v^*)^{\frac{1}{n}} =  v a^{\frac{1}{n}} v^*$.  Then
$$p_b = \lim_n \, b^{\frac{1}{n}} =
{\lim}_n \, v a^{\frac{1}{n}} v^* = v p_a v^* = v v^* ,$$
where these are weak* limits.
If $va \in A$ then $v B_p = v \overline{aBa}
\subset B$, where $p = p_a$.  Hence $B_p v^* \subset B$,
so that $a v^* \in B \cap A^{**} = A$.
\end{proof}

\begin{definition} \label{peddef} Let $a,b \in  {\mathfrak S}_A$.
We say that $a$ is {\em Pedersen equivalent} to $b$ in $A$,
and write $a \sim_A b$,  if the
equivalent conditions in the last theorem hold.
We will see some other characterizations of this variant of Pedersen equivalence
later in our paper.
Sometimes we write $a \sim_A b$ {\em via} $x$ and $y$, if
Theorem  \ref{varp} (ii)$^{\prime \prime}$  (or equivalently (ii)) holds.
Similarly if $n \in \Ndb$ we write $a^n  \sim_A b^n$
if (ii)$^{\prime \prime}$ in Theorem \ref{varp} holds with
$a$ and $b$ replaced by $a^n, b^n$.  We will see that
this happens iff $a \sim_A b$.  If $a, b \in
{\mathfrak c}_A$ with $\Vert a \Vert = \Vert b \Vert$,
then we say that $a \sim_A b$ if some positive scalar multiple of $a$ is Pedersen equivalent in the sense
above to the same multiple of $b$.
\end{definition}

It should be admitted at the outset that our Pedersen equivalence
(which generalizes the $C^*$-algebra variant--see Remark (2) below) happens
more rarely in general operator algebras than it does in
$C^*$-algebras.  However this is similar, for example,
 to the fact that
 right ideals with left contractive approximate identities
are rarer in general operator algebras than in $C^*$-algebras, but
are still interesting (they are related to important topics, 
and have a compelling theory--see e.g.\ \cite{BHN,BRead}.   While 
some algebras may have none, in other algebras they can be key for some purposes).
Results such as Proposition \ref{isu2},
 Theorem  \ref{anys}, and Theorem  \ref{isaq}
show that our Pedersen equivalence happens quite naturally.

\begin{xrem}
  (1) \   Of course if $a, b \in {\rm Ball}(A_+)$ in the case
 $A$ is a $C^*$-algebra, then
$a \sim_A b$ iff there exists $x \in A$ with $a = x x^* , b = x^* x$
(the usual Pedersen equivalence).   This is quite
obvious; for example if $a \sim_A b$, with notation as in
Theorem \ref{varp} (ii)$^{\prime \prime}$,
 then $y^* y = b^{\frac{1}{2}} b^{\frac{1}{2}} = b$
and $y y^* = a^{\frac{1}{2}} a^{\frac{1}{2}} = a$.

If $A$ is a $C^*$-algebra
then unitary equivalence of elements  $a,b \in
{\mathfrak S}_A$
implies $a \sim_A b$.  Indeed if $b = u a u^*$ for a unitary $u$ in
the multiplier algebra of $A$
(or simply a unitary in $A^{**}$ such that $au^*, ua \in A$),
then Ran$(b) = u {\rm Ran}(a)$.  Thus $u p_a = p_b u p_a$.
If we set $x = a^{\frac{1}{2}} u^*$ and $y =  u
a^{\frac{1}{2}}$ then
$x y =
a^{\frac{1}{2}} u^* u a^{\frac{1}{2}} = a$.
Also, $y x =  u a u^*  = b$, and
$x x^* =
a^{\frac{1}{2}} u^*u (a^{\frac{1}{2}})^*$.
So $a \sim_A b$.

(2) \ From Theorem \ref{varp} (iii) we see that in a finite dimensional operator algebra $A$,
$a \sim_A b$ iff $a \sim_r b$.  Here $a,b \in {\mathfrak S}_A$.
 We do not know if this is true in general.

(3) \ The relation $\sim_A$  in the case
that $A =  M_n$ is exactly unitary equivalence of elements
of ${\mathfrak S}_{M_n}$. The one direction
of this follows from (1).
For the other direction, if $a \sim_A b$ then
we can get $b = w^* a w$ for a tripotent $w$
that the reader can easily check may be replaced by a
unitary.  Unitary equivalence and similarity are the same
for positive matrices, but this is false for ${\mathfrak S}_{M_n}$
(even for $2 \times 2$ upper triangular matrices).

(4) \  We mention three other  interesting ways to state Pedersen equivalence: \\
(a) \ Theorem \ref{varp} (ii)$^{\prime \prime}$ is
equivalent with the same $x$ and $y$ to the same condition, but with $|x^{*}| = |c^{*}|$ replaced by
 $|x| = |d|$.   Clearly this is implied by (ii)$^{\prime \prime}$, and conversely it implies (ii),
 since $x = r(x) |x| = r(x) |d| = r(x) r(d)^* d$, which gives $x = x p_b$.

(b) \  $a \sim_A b$ iff there exist
$x,y \in A$ with $a = xy, b = yx$,
 and for each $n \in \Ndb$ there exist
contractions $x_n,y_n$ with $x = c^{1 -\frac{1}{n}} x_n, y = y_n c^{1-\frac{1}{n}}$.   The $x_{n}, y_{n}$ can be chosen
in $A$ if one wishes.  To see that
this is equivalent to (ii)$^{\prime \prime
 \prime}$, in the one direction simply  let $x_n = c^{\frac{1}{n}} R,
y_n = S c^{\frac{1}{n}}$.  For the other direction, if  subnets
$x_{n_k} \to R$ and $y_{n_k} \to S$ weak*, then (ii)$^{\prime \prime
 \prime}$  holds.  This is because
by the continuity of the disk algebra
functional calculus defining powers in \cite[Section 2]{BRead},
it is easy to argue that
$c^{1 - \frac{1}{n}} \to c$ in norm.

(c)  \ $a \sim_A b$ iff $a \sim_r b$ with some kind of `compatibility'
between the $x$ and $y$ elements for different roots.
That is, for
all $n \in \Ndb$ one may write $a^{\frac{1}{n}} = x_n y_n, b^{\frac{1}{n}} = y_n x_n$,
 where $x_n , y_n$ satisfy some conditions
 like $y_k x_{\ell} = y_r x_s$ whenever
$\frac{1}{k}
+ \frac{1}{\ell} = \frac{1}{r} + \frac{1}{s}$; or like $y_k x_{\ell} y_m = y_r x_s y_t$
whenever $\frac{1}{k}
+ \frac{1}{\ell} + \frac{1}{m} = \frac{1}{r} + \frac{1}{s}+ \frac{1}{t} \leq 1$.
To see for example that $a \sim_A b$ is implied by the last condition, together with
a couple of other reasonable conditions such as $x_n, y_n$ are contractions
and $y_n \in \overline{bAa}$:
Note  $y_{2n} x_{2n} y_2 = b^{\frac{1}{2n}} y_2$ converges to $y_2$, but it also
equals $y_n x_4 y_4 = y_n a^{\frac{1}{4}}$.  It follows that $(y_n a)$ converges,
so Theorem \ref{varp} (iii) holds.
Conversely, if Theorem \ref{varp} (ii) holds,
let $x_n = a^{\frac{1}{2n}} v^*, y_n = v a^{\frac{1}{2n}}$.
As in an argument in the proof of Theorem \ref{varp}, these are in $A$, $x_n y_n = a^{\frac{1}{n}},
y_n x_n = b^{\frac{1}{n}},$ and
$$y_r x_s y_t = v  a^{\frac{1}{2r}} a^{\frac{1}{2s}} v^* v
a^{\frac{1}{2t}} =  v a^{\frac{1}{2}(\frac{1}{r} + \frac{1}{s} + \frac{1}{t})} .$$
Thus  $y_k x_{\ell} y_m = y_r x_s y_t$ if  $\frac{1}{k}
+ \frac{1}{\ell} + \frac{1}{m} = \frac{1}{r} + \frac{1}{s}+ \frac{1}{t} \leq 1$.

(5) \ Note that $p_a$ and $p_b$ are the left and right
support projections of the element $x$ in Theorem \ref{varp} (ii) (and the
variants of (ii)).
\end{xrem}

The choice of the square root, as opposed to a different root,  in (i) and the variants of (ii) in Theorem \ref{varp} is
not essential;  one can instead take the $r$th power and the $(1-r)$th power, for any $r  \in (0,1)$.    For example:

\begin{proposition} \label{riok}  If
$a,b \in  {\mathfrak S}_A$ then $a  \sim_A b$  iff there exist
$x,y \in A$ with $a = xy, b = yx$ and
$|y| = |a^{r}|$ for some $r \in (0,1)$ and $x = x p_b$;   or iff there exist
$x,y \in A$ with $a = xy, b = yx$ and
$|y| = |a^{r}|$  and $|x^*| = |(a^{1-r})^*| $ for some $r \in (0,1)$.   \end{proposition} \begin{proof} Exactly as in the proof of Theorem \ref{varp}, but replacing
$c$ by $a^{r}$.   For example, given (iv) in Theorem \ref{varp},
let $y = va^r, x = a^{1-r} v^*$, then
it is easy to see, by the same arguments as in the proof of  Theorem \ref{varp}  that $x,y \in A$ with $a = xy, b = yx$, and
$y^* y = (a^r)^* v^* v a^r = (a^r)^* p_a a^r = |a^{r}|^2,$ so that  $|y| = |a^{r}|$, and $x = x p_b$.
Conversely, if
$x,y \in A$ with $a = xy, b = yx$ and
$|y| = |a^{r}|$, and $x = x p_b$, then all of the arguments in the first paragraph of the proof that
(ii) $\Rightarrow$ (ii)$^{\prime \prime}$ go through with $c$  replaced  by $a^{r}$.   In the second
  paragraph of that proof, one sees that $xva^r = a^{1-r} a^r$, hence
$xv a =  a^{1-r} a$, so that $xv = a^{1-r}$.  Also
$va = bv$, since $va = va^r  a^{1-r} = y  a^{1-r} = yx v = bv$.
  One can modify the lines that follow in the proof we are mimicking,
to get $y = v a^r = b^r v, p_b = v v^*,$ and $y y^* = a^r (a^r)^*$.  Thus
$a^{1-r} v^* = xv v^* = x$, and so
$x x^* = a^{1-r} (a^{1-r})^*$.  This all is more than enough
to obtain  (iv) in Theorem \ref{varp}, indeed $va = va^r a^{1-r} = y a^{1-r} \in A,
a v^* = a^r a^{1-r} v^* = a^r x \in A$.  Since $v a = bv,$ we have $b = v a v^*$.
\end{proof}

\begin{proposition} \label{ifr}  Suppose that
$a,b \in  {\mathfrak S}_A$ and $r \in (0,1) \cup \Ndb$.
Then $a  \sim_A b$ iff $a^r \sim_A b^r$.
\end{proposition}

\begin{proof}  If $a \sim_A b$, then using notation
from the proof of Theorem \ref{varp} above, let $y' = v c^r,
x' = c^r v^*$.  Note that $y' = v c^r \in v \overline{a A a} \subset A$
since $y = v c \in A$.  Similarly $x' \in A$.
Then $x' y' = c^{2r} = a^r$, and $y' x' =
v c^{2r} v^* = v a^r v^* = b^r$.  The latter follows
by Lemma \ref{vpow} and
 Lemma \ref{cofr}.  So $a^r \sim_A b^r$.

The converse is very similar, suppose that
 $a^r \sim_A b^r$ via $a^{\frac{r}{2}} v^*,
v a^{\frac{r}{2}}$, so that $b^r = v a^r v^*$.  Let $x = a^{\frac{1}{2}} v^*, y =
v a^{\frac{1}{2}}$.  Then $x,y \in A$ and
$x y = a$ as before, and $$yx = v a v^* =
 ((v a v^*)^r)^{\frac{1}{r}} = (v a^r v^*)^{\frac{1}{r}}
= (b^r)^{\frac{1}{r}} = b,$$ as can be
easily argued following the idea in the last paragraph.
\end{proof}

\begin{corollary} \label{ieq}  $\sim_A$ is an equivalence relation
on ${\mathfrak S}_A$.
\end{corollary}

\begin{proof} We leave this as an exercise using e.g.\ Theorem \ref{varp} (iv).
 \end{proof}

\begin{theorem}  \label{vpr}  Suppose that
$a \in  {\mathfrak S}_A$, and let
$c = a^{\frac{1}{2}}$.  If there exists
$x,y \in A$ with $a = xy$, and
$|y| = |c|$ and $|x^*| = |c^*|$, then
$b = yx$ is in ${\mathfrak S}_A$,
and $a \sim_A b$.
\end{theorem}

\begin{proof}   Using the notation and proof of  (ii)$^{\prime}$
 $\Rightarrow$ (ii) in
Theorem \ref{varp} above,
we have
$y = vc$,
and hence $v v^* b = b$,
and
$x = cw$
and $b w^* w = b$.  Following that proof, we see that
$w = v^*$ again, so that
$b = v v^* b v v^*$.  So
$$\Vert 1 - 2 b \Vert =
\Vert v v^* - 2  v v^* b v v^* \Vert \leq \Vert 1 - 2  v^* yx v  \Vert =
\Vert  1 - 2 c^2 \Vert =
\Vert  1 - 2 a \Vert.$$
Similarly, $\Vert 1 -  b \Vert \leq \Vert 1 - a \Vert$.
So $b = yx$ is in ${\mathfrak S}_A$,
and $a \sim_A b$.  \end{proof}

\section{Open tripotents for operator algebras and
Peligrad-Zsid\'o equivalence} \label{optoo}

The partial isometry or tripotent
$v$ that occurs in Theorem \ref{varp} is central
to our study of the above  variant of Pedersen equivalence.
It also happens to be what we will call a $*$-{\em open tripotent}.
 In the present section we define this notion and develop the associated theory.
 Later we will apply this theory back to Pedersen equivalence,
 and to coarser notions of equivalence and subequivalence.

In the introduction we discussed the notion from \cite{BN0}
of an open tripotent in $B^{**}$, and gave several equivalent definitions
of these objects.  We will also use below the following
characterization, which is proved below in the  Appendix (Proposition \ref{pzis}):
A partial isometry $v$  in $B^{**}$  is an open tripotent in the
sense of \cite{BN0} iff
$p = v^* v$ and $q = v v^*$ are open projections in the sense of Akemann,
and $v B_p \subset B$ and $v^* B_q \subset B$.  As we said in the introduction,
first, $p$ and $q$ are said to be Peligrad-Zsid\'o equivalent in
this case, and
second, the conditions involving $q$ in the last line are automatic and are stated
for the sake of symmetry.

We define a {\em $*$-open tripotent} (or {\em $*$-open partial isometry})
{\em in} $A^{**}$, for
an operator algebra $A$, to be an open tripotent (in the sense of
\cite{BN0}) in $B^{**}$, which lies in $\Delta(A^{\perp \perp})$.
Here $B$ is a $C^*$-algebra generated by $A$.

\begin{lemma} \label{bn0l}  If $v$ is an open tripotent in
$B^{**}$ then so is $v^*$.  Hence if $v$ is a $*$-open tripotent in
$A^{**}$ then so is $v^*$.
\end{lemma}  \begin{proof}  We leave this as an exercise.
\end{proof}

Being a
$*$-open tripotent  makes sense
independently of the particular $C^*$-algebra $B$ generated by $A$,
since one may define `$*$-openness' in $A^{**}$ purely in terms
of $A$:  $u$ is a $*$-open tripotent in $A^{**}$
iff $u$ is a tripotent in the $C^*$-algebra
$\Delta(A^{**})$, and the projection
$$\hat{u} = \frac{1}{2} \left[ \begin{array}{ccl} u u^* & u \\
u^* & u^* u \end{array} \right]$$
in $M_2(\Delta(A^{**}))$, is an open projection
in $M_2(A)^{**}$
in the sense of \cite{BHN} (that is, there exists a net
$a_t \in M_2(A)$ such that $a_t = \hat{u} a_t = a_t \hat{u}
\to \hat{u}$ weak* in $M_2(A)^{**}$).  This is also equivalent to
$\hat{v}$ being open
in $M_2(A)^{**}$
in the sense of \cite{BHN}, where $v = u^*$.

  \begin{xrem} (1) \     Open projections in $A^{**}$ in the sense of
\cite{BHN}, are clearly $*$-open tripotents in $A^{**}$.  This
follows from the definition above since open projections in $B^{**}$
are open tripotents in the sense of \cite{BN0}.

(2) \ It is easy to see how to define a $*$-open tripotent in the second dual of a
strong Morita equivalence bimodule $X$ in the sense of \cite{BMP}, and a matching Peligrad-Zsid\'o equivalence for open projections in the second duals of the algebras acting on the left and right of $X$.  Most of the results below
about $*$-open tripotents and Peligrad-Zsid\'o equivalence will be almost identical in this more general setting, with almost identical proofs, and we leave this to the reader.  This will then simultaneously generalize much of the theory in the present paper, and the
theory of open tripotents in \cite{BN0}.   We have chosen not to present this here,
since it is often better to present a more concrete theory than
to present the most general abstract formulation.  We do remark however that many results in this more general
theory are derivable from the case we consider by taking
$A = {\mathcal L}(X)$, the linking algebra of $X$. \end{xrem}

We say that two projections $p, q \in A^{**}$ are
{\em Peligrad-Zsid\'o
equivalent in} $A^{**}$, and write
$p \sim_{PZ,A} q$ if there exists an open tripotent
 (in the sense of \cite{BN0})
    $v \in B^{**}$, which lies in $\Delta(A^{\perp \perp})$ (hence is $*$-open),
 implementing a Peligrad-Zsid\'o
equivalence of $p$ and $q$ in $B^{**}$.
That is, $v^* v = p, v v^* = q$ (since
$v$ is open this forces $p, q$ to be open projections
\cite{BN0}).

    \begin{xrem}  If $p \sim_{B,PZ} q$ it does not follow that
 $p \sim_{PZ,A} q$ necessarily.
   A simple counterexample is the diagonal matrix units
in the case $A$ is the upper triangular 2 by 2 matrices.
\end{xrem}

To connect our definition to  Pelegrad and Zsid\'o's original equivalence,
 as we did in the selfadjoint case (see Proposition \ref{pzis}), we have:

\begin{lemma} \label{chpz}  Two
projections $p, q \in A^{**}$ are
 Peligrad-Zsid\'o equivalent in $A^{**}$ iff
$p, q$ are open projections and there exists
$v \in B^{**}$ such that
$v^* v = p, v v^* = q$, and $vA_p \subset A$ and $v^* A_q \subset A$.
This is also equivalent to:
$p$ is open and there exists
$v \in \Delta(A^{**})$ such that
$v^* v = p, v v^* = q$, and $vA_p \subset A$.
\end{lemma}

\begin{proof}  If $p \sim_{PZ,A} q$ via $v$ then
as we said above, $p, q$ are open.
Also by Proposition \ref{pzis} from the Appendix,
$vA_p \subset vB_p \subset B \cap A^{\perp \perp}
= A,$ and similarly $v^* A_q \subset A$.

Conversely, suppose that $p, q$ are open and
that $vA_p \subset A$ and $v^* A_q \subset A$.  Then $v e_t
\in A$ if $(e_t)$ is a cai for
$A_p$.  Taking a weak* limit, $v p = v \in A^{\perp \perp}$.
Similarly, $v^* \in A^{\perp \perp}$, so
$v \in \Delta(A^{\perp \perp})$.
Now $(e_t)$ is also a cai for $B_p$ by \cite[Theorem 2.9]{BHN},
so it follows that for any $c \in B_p$ we have
$vc = \lim_t v e_t c \in AC \subset B$.  That is,
$vB_p \subset B$.  Similarly, $v^* B_q \subset B$, so
by Proposition \ref{pzis} from the Appendix,  $v$ implements a Peligrad-Zsid\'o
equivalence of $p_a$ and $p_b$ in $B^{**}$.
If $p$ is open and $vA_p \subset A$ then the argument above gives
$vB_p \subset B$.  So by  \cite[Lemma 1.3]{PZ}, $q$ is open and
$v^* B_q \subset B$.  Thus if $v^* \in A^{\perp \perp}$ then
$v^* A_q \subset B \cap A^{\perp \perp} = A$.
\end{proof}

\begin{xrem}
     Under the equivalent conditions of the
Lemma, we also have
 $A_q v \subset B_q v = (v^* B_q)^* \subset B$,
so $A_q v \subset B \cap A^{\perp \perp} = A$.
Similarly, $A_p v^* \subset A$.
\end{xrem}

\begin{theorem} \label{mat1}  Let $u$ be a tripotent in $A^{\perp \perp}$,
and let $$A_u = \{ a \in A : a = u u^* a u^* u \} .$$  Suppose that $u$
 is in the weak* closure of $A_u$, and that $A_u u^* \subset A$.  Then $u$
is an open tripotent in $B^{**}$.  If, further, $u^* \in A^{\perp \perp}$
then $u$ is a $*$-open tripotent in $A^{**}$, $u^* A_u \subset A$, and the weak* closure of
$A_u$ equals $$\{ \eta \in A^{**} : \eta = u u^* \eta u^* u \} =
u u^* A^{**} u^* u .$$

Conversely, if $u$ is a $*$-open tripotent in $A^{**}$ then
$u$ is in the weak* closure of $A_u$, and  $A_u u^* \subset A$.
\end{theorem}  \begin{proof}  We have $u A_u^* A_u \subset B$.  As explained in the introduction, $B^{**}_2(u)$ is a
$C^*$-algebra with product $x u^* y$ and involution $u x^* u$.
In this product and involution the $C^*$-algebra generated by
$u A_u^* A_u$ is contained in $B$.  The increasing cai for
this $C^*$-algebra converges to  a projection in
the $C^*$-algebra  $B^{**}_2(u)$, hence to a tripotent $w \in B^{**}$.
Since $w u^* w = w$, we have $w \leq u$ in the ordering of
tripotents.  It is an exercise to check that $w$ satisfies the
definition from \cite{BN0} of an open tripotent in $B^{**}$.
If $x, y \in A_u$ we have $u x^* y = w u^* u x^* y = w w^* u x^* y = w x^* y$.
Replacing $x$ with $x_t  \in A_u$ such that $x_t \to u$ weak*, we have $u u^* y =
w u^* y = w w^* y$.  Replacing $y$ with $x_t$ in the limit we have
$u = u u^* u = w w^* u = w$.  So $u$ is open in $B^{**}$.

If further $u^* \in A^{\perp \perp}$ then $u^* A_u \subset u^* B_u
 \subset B \cap A^{**} =A$.  Clearly $u$ is $*$-open  in $A^{**}$.
The weak* closure of
$A_u$ is contained in $u u^* A^{**} u^* u$ clearly.  Conversely, if $p = u^* u$,
which is open in $A^{**}$ by \cite{BN0},  and
$\eta \in u u^* A^{**} u^* u$, then $u^* \eta \in p A^{**} p$, so
there is a net $(a_t)$ in $_pA_p$  converging weak* to $u^* \eta$.
Then $u a_t \to  \eta$ weak*, and $u a_t \in  A$, indeed $u a_t \in  A_u$.  So the weak* closure of
$A_u$ is $u u^* A^{**} u^* u$.

For the converse, if $u$ is a $*$-open tripotent in $A^{**}$
then by the fact above the theorem there exists a
net $x_t \in A_u$ with $x_t \to u$ weak*.  If $x \in A_u$ then $x \in B_u$.
In the notation of \cite{BN0}, $B_u u^* = B(u) u^* \subset B$ by
facts in \cite{BN0}.  Then
$x u^* \in B_u u^* \cap A^{**} \subset B \cap A^{**} = A$.
  \end{proof}

In the following, ${\mathfrak S}_{A_u}$ is computed with
respect to operator algebra $A_u$ with the Peirce product.  That is, ${\mathfrak S}_{A_u} = \{ a \in A_u :
\Vert u - K a \Vert \leq 1 \}$, where $K = 1$ or $2$.

\begin{lemma} \label{ord}  If $u$ is
a $*$-open tripotent in $A^{**}$ then $A_u$ is a subalgebra of the
$C^*$-algebra $B(u)$ with its `Peirce product', and it has a cai
with respect to this product.   If $u, w$ are
$*$-open tripotents in $A^{**}$, then TFAE:
\begin{enumerate} \item [(i)] $u \leq w$ in the ordering of tripotents,
\item [(ii)]   $A_u$ is a subalgebra (indeed is a HSA) of $A_w$ in the product above,
\item [(iii)]  ${\mathfrak S}_{A_u} \subset {\mathfrak S}_{A_w}$  as sets,
\item [(iv)]  $\hat{u} \leq \hat{w}$.
\end{enumerate}
Indeed a tripotent in $A^{**}$ which
 is  dominated by $w$  in the ordering of tripotents, is $*$-open in $A^{**}$ iff it is
an open projection in the operator algebra $A_w^{**}$ (with its
Peirce product).
Moreover, $u = w$ iff $A_u = A_w$ as sets and as algebras
(with the product above).    \end{lemma}  \begin{proof}
That $A_u$ is a subalgebra  of $B(u)$
 follows from the result $A_u u^* \subset A$ above.
Since the second dual of $A_u$ is its weak* closure
$u u^* A^{**} u^* u$, which contains $u$ which is an identity
in the `Peirce product', it follows from e.g.\ \cite[Proposition 2.5.8]{BLM}
that $A_u$ has a cai.

 If $u \leq w$ then
it is clear that $A_u \subset A_w$ as sets.

(i) $\Rightarrow$ (ii) \ If $x, y \in A_u$
then $x w^* y = x u^* u w^* y = x u^* y$, so
$A_u$ is a subalgebra of $A_w$.  We leave it as an exercise that it is a HSA.

(ii) $\Rightarrow$ (i) \ If
$A_u$ is a subalgebra of $A_w$ then $(A_u)^{**}$ is a
subalgebra of $(A_w)^{**}$.  Thus $u w^* u = u$, so $u \leq w$.

(i) $\Leftrightarrow$ (iv) \ is in \cite[Proposition 2.7]{BN0}.

(ii) $\Rightarrow$ (iii) \  If (ii) holds then by (i) $u$ may be viewed as a projection
in $(A_w)^{**}$ dominated by the identity $w$.  Let $x \in A_{u} \subseteq A_{w}$ with  $||u-2x|| \leq 1$.  Then
$||w-2x|| =||w-u+u-2x|| \leq \max \{1,||u-2x|| \} \leq 1$.

(iii) $\Rightarrow$ (i) \ Taking second duals, and using
\cite[Lemma 8.1]{BRead} we have ${\mathfrak S}_{A^{**}_u} \subset {\mathfrak S}_{A^{**}_w}$.  By  the fact mentioned in the introduction  that ${\mathfrak c}_A$ densely spans $A$,
we have that $A^{**}_u \subset A^{**}_w$.
In particular,  $\Vert w - 2u \Vert
\leq 1$ since $\Vert u - 2u \Vert
\leq 1$.
We wish to show that $u w^* u = u$, or equivalently that $u$
 is a projection in the $C^*$-algebra $C = B^{**}_2(w)$.  Working inside the latter
 $C^*$-algebra we have $w = 1$.  Thus it suffices to prove
 that any partial isometry $V$ in a unital $C^*$-algebra
 $C$, with $\Vert 1 - 2 V \Vert \leq 1$, is a
 projection.  Since $V \in {\mathfrak S}_{C}$, its left
 support and right supports agree, so $V^* V = V V^*$ equals  a
 projection $p = s(V)$.  If $D = pCp$ then
 $V$ is a unitary in $D$ with $\Vert p - 2V \Vert \leq 1$.  By spectral theory it follows that $V = p$.  This
 completes the proof of (i).

Given an open projection $p \in (A_w)^{\perp \perp}$,
we recall that $A_w$ is a subalgebra of $B(w) \subset
B^{**}_2(w)$.  So by \cite[Theorem 2.4]{BHN},
$p$ is open in $B^{**}_2(w)$.  By \cite[Corollary 2.11]{BN0}, $p$
is an open tripotent $u$ in $B^{**}$.  We also have
that $u^* = p^* = w^* p w^* \in A^{**}$.  So $u$
is a $*$-open tripotent in $A^{**}$.

Conversely, any $*$-open tripotent $u \in A^{**}$  is open as a tripotent in $B^{**}$.
  If $u \leq w$, then $u \in (A_w)^{\perp \perp} \subset B^{**}_2(w)$.
By \cite[Corollary 2.11]{BN0}, $u$ is open as
a projection in $B^{**}_2(w)$.  Since
$A_w$ is a subalgebra of $B(w)$ and $u
 \in (A_w)^{\perp \perp}$,
$u$ is open as
a projection in $(A_w)^{\perp \perp}$ by
\cite[Theorem 2.4]{BHN}.   \end{proof}

 \begin{xrem}  If $A_u = A_w$ as subsets of $A$, it need not follow
that $u = w$.   This may be seen by considering for example any two distinct
unitaries in $K(H)^{**} = B(H)$.   Indeed $A_u = A_w$ iff $u^* u = w^* w$
and $u u^* = w w^*$.  To see the harder direction of this, if
$A_u = A_w$ then in the second duals we see that $u \in B^{**}_2(w)$,
so that $u^* u \leq w^* w$ and $u u^* \leq w w^*$.  Similarly, $w^* w \leq u^* u$
and $w w^* \leq u u^*$. \end{xrem}

From Sections 2 and 3 of \cite{BN0} one can write down a host of
properties and behaviors of open tripotents in  $A^{**}$.  For example:

\begin{lemma} \label{sups}  An increasing net of $*$-open tripotents
in $A^{**}$ has a least upper bound tripotent in $A^{**}$, namely
its weak* limit, and this limit tripotent is  $*$-open in $A^{**}$.
\end{lemma}

\begin{proof}  Clear from \cite[Proposition 2.12]{BN0} and the
characterization of $*$-open tripotents
in $A^{**}$ as open tripotents in $B^{**}$ which are also in the
von Neumann algebra $\Delta(A^{**})$.
 \end{proof}

\begin{lemma} \label{sups2}   If $u$ and $v$ are
$*$-open tripotents in $A^{**}$, which are both
dominated by some tripotent in $B^{**}$,
then there exists a smallest  tripotent $u \vee v$ in $B^{**}$
which dominates both $u$ and $v$, and this
tripotent is a $*$-open tripotent in $A^{**}$.
\end{lemma}

\begin{proof}  Of course $u$ and $v$ are open tripotents in $B^{**}$.
That there exists a smallest  tripotent $u \vee v$ in $B^{**}$
which dominates both $u$ and $v$ is well known (see e.g.\
\cite[Lemma 3.7]{BN0}).   By \cite[Lemma 3.7]{BN0}, $\widehat{u \vee v}
= \hat{u} \vee \hat{v}$.   However $\hat{u}$ and  $\hat{v}$ are open
projections in $M_2(\Delta(A^{**}))$, so we deduce
 that  $\widehat{u \vee v}$ is an open projection
and is in $M_2(\Delta(A^{**}))$, so that $u \vee v$ is
a $*$-open tripotent in $A^{**}$.
\end{proof}

\begin{lemma} \label{sups3}   A family $\{ u_i \}$ of $*$-open tripotents
in $A^{**}$, which are all dominated by a tripotent in $B^{**}$,
has a least upper bound amongst the tripotents in $B^{**}$,
and this  is a $*$-open tripotent  in $A^{**}$.
\end{lemma}

\begin{proof} Of course the $u_i$ are open tripotents in $B^{**}$,
so by \cite[Proposition 3.8]{BN0} the family has a least upper bound $v$
amongst the
tripotents in $B^{**}$,
and this  is an open tripotent  in $B^{**}$.  The supremum of any
finite subfamily of  $\{ u_i \}$ is a
$*$-open tripotent by Lemma \ref{sups2}, and such supremums
form an increasing net with least upper bound $v$.
so $v$ is a $*$-open tripotent  in $A^{**}$ by Lemma \ref{sups}.
\end{proof}

We now rephrase Pedersen equivalence in terms of $*$-open
tripotents.   We can rephrase a special case of Lemma
\ref{chpz} as follows.  Let $v \in B^{**}$ with $v^* v = p_a$ for some
$a \in {\mathfrak S}_A$ with $va \in A$.  Then
$v B_{p_a} = v \overline{aBa} \subset B$,
and so $v$ is an open tripotent in $B^{**}$, and $v = \lim_n v a^{\frac{1}{n}} \in A^{**}$.
In this case, $v$ is a $*$-open tripotent in $A^{**}$ iff $v^* \in A^{**}$, and iff
$a v^* \in A$.

\begin{proposition} \label{simt}If $a,b \in {\mathfrak S}_A$ then  $a \sim_{A} b$ iff
there exists a
$*$-open tripotent $v$ with $v^* v = p_a, v v^* = p_b$ and $b = vav^{*}$.
Also, in this case the $x$ and $y$ in Theorem {\rm \ref{varp} (ii)} are in bijective
correspondence with such
$*$-open tripotents $v$
satisfying all the conditions in the last line.    The bijection
here takes $v$ to $y = vc$ (and $x=cv^{*}$), and its inverse takes $y$ to $r(y)
r(c)^{*}$.  Here $c = a^{\frac{1}{2}}$.  \end{proposition}

\begin{proof} The first assertion is clear from Theorem \ref{varp} (iv) and Lemma
\ref{chpz} (and the fact that $v a \in A$ implies that $v \, \overline{aAa} \subset A$).
  To finish we show that the
two maps described are inverses of each other.  If $y$ is as in Theorem {\rm
\ref{varp} (ii)} we have  $r(y) r(c)^{*} c = r(y) |c| = r(y) |y| = y$.   Conversely,
if $v$ is a $*$-open tripotent satisfying all the conditions in the first line of
the theorem statement, then $vc = v r(c) |c| = v r(c) |vc|$.   Now $v r(c)$ is a
tripotent, with left support $v s(c)v^{*} = vv^{*}$ and
right support $r(c)^{*} v^{*} v r(c)= p_{a}$.  By the `uniqueness of the polar decomposition'
 we deduce that
$v r(c)=r(vc)$.  Thus $r(vc)r(c)^{*} = v r(c)^{*} r(c) = v p_{a} = v$.  This
completes the proof. \end{proof}

\begin{proposition} \label{isu2}  Let
 $v$ be  a $*$-open tripotent in $A^{**}$,
and suppose that $p = v^* v$ is the support projection for some
$a \in {\mathfrak S}_A$.
Then $q = v v^*$ is the support projection for some $b \in {\mathfrak S}_A$,
and $a \sim_A b$.
 \end{proposition} \begin{proof}    Let $b = v a v^*$.
As in \cite{BRead},  $\overline{aAa} =  A_p$ since $p = s(a)$.
Note that
$b = v a^{\frac{1}{2}} a^{\frac{1}{2}}
v^* \in v A_p A_p v^* \subset A_q$, using the
Remark after Lemma \ref{chpz}.  Also,
$b \in {\mathfrak S}_A$ and
 $(v a v^*)^{\frac{1}{n}} =  v a^{\frac{1}{n}} v^*$
by Lemma \ref{vpow}.  As at the end of the proof of
Theorem \ref{varp},
$p_b =
v v^* = q$.
That $a \sim_{A} b$ now follows from
Proposition  \ref{simt}.        \end{proof}

\begin{xrem} A similar result: suppose that
$a \in  {\mathfrak S}_A$, and that for all $n \in \Ndb$, there exist
$x_n,y_n \in {\rm Ball}(A)$ with $a^{\frac{1}{n}} = x_n y_n$,
and  the sequence $( y_n a)$ has a norm convergent subsequence
with limit $va$ say, where $v$ is a weak* cluster point of $(y_n)$.
Then $v$ is a $*$-open tripotent, $v \in \Delta(A^{**})$, and $b = v a v^* \in {\mathfrak S}_A$ with
$a \sim_A b$.
To prove this, suppose that $x_{n_k} \to S$ and $y_{n_k} \to v$ weak*, and $y_{n_k} a \to va$ in norm.
As usual we may assume that $y_n = y_n p_a, x_n = p_a x_n$, and $S = p_a S$ and $v = v p_a$.
Then  $x_{n_k} y_{n_k} a \to S v a = p_a a = a$, so that $Sv = Sv p_a = p_a$.
This forces $S^* = v$, and $v$ is a tripotent in $\Delta(A^{**})$.
Since $y_{n_k} a \to va$ in norm, $va \in A$.  As in Lemma
\ref{chpz} above, $vB_p \subset B$,
where $p = p_a$.  By \cite[Lemma 1.3]{PZ}, if $q = v v^*$ then
$q$ is open and $v^* B_q \subset B$.   So $v$ is open, and
$a^{\frac{1}{2}} v^* \in B_p v^* \cap A^{\perp \perp} = A$,
and similarly $v a^{\frac{1}{2}} \in A$.   So by Proposition
\ref{isu2} we have $b = v a v^* \in {\mathfrak S}_A$ and $a \sim_A b$.
\end{xrem}

\section{Hereditary bimodules} \label{hbm}

Hereditary bimodules, defined below, as well as being
interesting in their own right,  will be  important for us for two main reasons.
First,  they are useful in understanding Pedersen and Blackadar equivalence: if one examines the
algebraic structure of $\overline{aAb}$ when $a$ and $b$ are
equivalent in these ways, then one is led naturally to hereditary bimodules.    (Indeed we give several characterizations of these equivalences in terms of hereditary bimodules in this section and the next.)
Secondly, hereditary bimodules will turn out to be much the same thing
as $*$-open tripotents $u$, via the correspondence
$u \mapsto A_u$, with $A_u$ defined as in the previous section.
One important definition of open tripotents from \cite{BN0} is as
the `support tripotents' of what was called an {\em inner $C^*$-ideal}.
Hereditary bimodules are the analogue of inner $C^*$-ideals in our setting,
and we shall see in Proposition  \ref{mortri} below
that $*$-open tripotents are the same as `support tripotents'
of hereditary bimodules.

We define a {\em hereditary bimodule} to be
the $1$-$2$ corner of a HSA in $M_2(A)$, where that
HSA is isomorphic to $M_2(C)$ for an
approximately unital operator  algebra $C$,
via a `corner preserving'
 completely isometric homomorphism $\rho : M_2(C) \to
M_2(A)$.   Thus we have subalgebras $D, E$ of $A$,
and subspaces $X, Y$ of $A$, such that
the subalgebra ${\mathcal L}$ of $M_2(A)$ with rows $C,X$ and $Y,D$, is
a HSA in $M_2(A)$, and there is a
completely isometric isomorphism $\rho$ from $M_2(C)$ onto ${\mathcal L}$
taking each of the corners $D, E, X, Y$ of ${\mathcal L}$ onto
the copy of $C$ in the matching corner of $M_2(C)$.
In this case we say that $X$ is a
{\em hereditary bimodule},
or a {\em hereditary} $D$-$E$-{\em bimodule},
and we call ${\mathcal L}$ the
{\em linking algebra of the hereditary bimodule}.  It then
follows easily that $(D,E,X,Y)$ is a {\em strong Morita equivalence context}
 in the sense of \cite{BMP}, that $Y$ is the `dual module'
$\tilde{X}$ considered in that theory, $X = \tilde{Y}$, and that ${\mathcal L}$ is the
 linking algebra in the sense of \cite{BMP} for this strong Morita equivalence.
Conversely, given a
strong Morita equivalence context $(D,E,X,Y)$ whose
linking algebra ${\mathcal L}$ is a HSA in $M_2(A)$,
then saying that ${\mathcal L}$ is isomorphic to $M_2(C)$
in the way described above,
 is, by \cite[Corollary 10.4]{Bmuls},  the same as saying that the $1$-$2$ corner $X$ is
linearly completely isometric to an
approximately unital operator  algebra.
That is, given a
strong Morita equivalence context $(D,E,X,Y)$ whose
linking algebra ${\mathcal L}$ is a HSA in $M_2(A)$,
then $X$ is a hereditary bimodule
 iff $X$ is
linearly completely isometric to an approximately unital operator  algebra.

We point out that $D = XY$ and $E = YX$ above.  These are
HSA's in $A$.  But also
$X = DAE$ and $Y = EAD$, and these are inner ideals in $A$ (that is,
$X A X \subset X$ and similarly for $Y$).  So $X$ and $Y$ are retrievable from
$D$ and $E$.  To see this, note that clearly $X = D X E \subset DAE$, and conversely
$DAE = XY A YX \subset X$ since $X$ is an inner ideal.   So
$X = DAE$, and similarly $Y = EAD$.

In the notation above, note that a hereditary bimodule $X$ is
an approximately unital operator algebra in the product transferred
from $C$ via $\rho$.  The weak* limit in $A^{**}$ of the cai for
$X$ will be called a {\em support tripotent} of $X$ (in the proof
of the next result it will be seen to be a tripotent).  This support tripotent is
not unique, as we discussed for example in  the remark before Lemma \ref{sups}, but
can be easily be made so similarly to the case in \cite{BN0} by keeping track
of one further peice of structure on a hereditary bimodule, such as its `product'
(or one of the other items in Lemma \ref{ord}).

We define a {\em principal hereditary  $A_a$-$A_b$-bimodule},
or {\em principal hereditary bimodule} for short,  to be
a hereditary  $A_a$-$A_b$-bimodule $X$ for some
$a,b \in {\mathfrak S}_A$.
We have $X = \overline{aAb}$ by Proposition \ref{aab}.
The
linking algebra ${\mathcal L}$ of the bimodule is the subalgebra of $M_2(A)$ with
rows $\overline{aAa},\overline{aAb}$ and $\overline{bAa}, \overline{bAb}$.
Note that this ${\mathcal L}$  is certainly a HSA in $M_2(A)$.

\begin{xrem}
  For any $a,b \in {\mathfrak S}_A$,
the subalgebra of $M_2(A)$ with rows $\overline{aAa},\overline{aAb}$ and $\overline{bAa}, \overline{bAb}$,
is a HSA in $M_2(A)$ clearly.  However it will not in general correspond to a strong Morita equivalence.
Indeed, if $X = \overline{aAb}$, there may be no HSA in $M_2(A)$ with $X$ as its $1$-$2$ corner
which  corresponds to a strong Morita equivalence.  Obvious counterexamples exist
in the case $A$ is the upper triangular $2$ by $2$ matrices.  \end{xrem}

\begin{proposition}  \label{mortri}  If $u$ is
a $*$-open tripotent in $A^{**}$, then $A_u$ is a hereditary bimodule.
Conversely, any hereditary  bimodule
equals $A_u$ for a $*$-open tripotent $u$; and $u$ is
a support tripotent of this hereditary  bimodule.  \end{proposition}

\begin{proof}  If $u$ is a $*$-open tripotent in $A^{**}$ then by
Lemma \ref{bn0l} so is $v = u^*$.  Let $p = v^* v, q = v v^*$.
The subalgebra ${\mathcal L}$ of $M_2(A)$ with
rows $_pA_p, \, _pA_q$ and $_qA_p, \, _qA_q$,
is a HSA of $M_2(A)$.  Note $_pA_q = A_u$.   Since
$A_u$ is an approximately unital
operator algebra in the `Peirce product', it is a
hereditary bimodule by the definition
of the latter.   Explicitly, if  $S = {\rm diag}(1,v)$, then
$\ell \mapsto S^* \ell S$ is a completely contractive homomorphism
from ${\mathcal L}$   to $M_2(_pA_p)$, whose
completely contractive inverse is $k \mapsto S k S^*$.

For the converse, we use the notation in our definition
and discussion of `hereditary bimodule' at the start of
Section \ref{hbm}.
Let $e$ be the identity of $C^{**}$.
Note that $\rho^{**} : M_2(C^{**}) \to M_2(A^{**})$ is a
completely isometric corner preserving homorphism too, and
it restricts to an injective $*$-homomorphism from
$M_2(\Cdb e)$ into $\Delta({\mathcal L}^{**})$.
Note that the identity of ${\mathcal L}^{**}$ is
$p + q$ where $p$ is the support projection
of $D$ and $q$ is the support projection of $E$.  Also
$p$ and $q$ correspond via $\rho^{**}$ to $e \oplus 0$ and $0 \oplus e$
in $M_2(C^{**})$.
Let $u$ be the image in $A^{**}$ of $e$ under the $1$-$2$ corner of
$\rho^{**}$.  Then $v = u^*$ is the image in $A^{**}$ of $e$ under the $2$-$1$ corner
 of $\rho^{**}$.  Also,
 $u$ is a tripotent in $X^{\perp \perp}$, and $v v^* = q, v^* v = p$
(since the matching facts are true in $M_2(\Cdb e)$).
Also, $u A_p \subset A$, since the matching
products are in $C$, and $\rho(C) \subset M_2(A)$.
So $u$ is a $*$-open tripotent in $A^{**}$  by Lemma \ref{chpz}
(one could also use Theorem \ref{mat1} to see this).
Since a cai for $C$ converges weak* to $e$, it is clear that
$u$ is a support tripotent in the sense above.
  \end{proof}

The following is an intrinsic characterization of hereditary bimodules:

\begin{theorem} \label{phii}   An inner ideal $X$ in $A$ is
a hereditary bimodule, or equivalently is of the
form $A_u$ for a $*$-open tripotent $u$, if and only if there
are nets of contractions $c_n \in X$ and $d_n \in A$, such that
\begin{enumerate} \item [{\rm (1)}]  $c_n d_n x \to x$ and $x d_n c_n \to x$
in norm for all $x \in X$,
\item [{\rm (2)}]  $(d_n c_m)$ and $(c_m d_n)$ are norm convergent
with $n$ for fixed $m$, and are norm convergent with $m$ for fixed
$n$.
\end{enumerate}  If these hold then a subnet of the $c_n$ converge weak* to a
support tripotent $u$ for $X$, so that $X = A_u$.
 \end{theorem}

\begin{proof}  ($\Rightarrow$) \  Suppose that $u$ is a  $*$-open tripotent, so that $A_u$
 is a hereditary bimodule, isomorphic to an approximately unital operator algebra
 $C$.  Then since these conditions clearly hold in $C$, by
applying the isomorphism $\rho$ in the definition of `hereditary bimodule' we see that they will hold too in $M_2(A)$.

 ($\Leftarrow$) \ Replacing by subnets, we may suppose that $c_n,
 d_n, c_n d_n,$ and $d_n c_n$, converge in the weak* topology
to say elements $u, w, p, q$, respectively.
Then $c_n d_n c_m$ converges  weak* to $u wc_m = c_m = pc_m$.
Taking the limit over $m$ we see that $u wu = u = p u$, so
$(uw)^2 = u w  = pu w$.  Since $u w$ is a contraction, it is a
projection.  Also, $c_n d_n c_m d_m$ converges  weak* to $u wc_m d_m
= c_m d_m = p c_m d_m$.  Taking the limit over $m$ we see that
$u w p = p = p^2$.   So $p$ is a projection too, and now we see that
$p = uw$.   Also $u = pu$, and similarly $wp = \lim_n d_n p = w$.
It follows that $u = w^*$, and this is a partial isometry, and
is in $\Delta(A^{**})$.  Similarly, $q = u^* u$.
Condition (1) implies in the limit that
 $p x = x = xq$ for $x \in X$, so $X \subset \, _pA_q$.
Hence $u$ is in the weak* closure of $A_u$.
Conversely, if $y \in \, _pA_q$ then
$u^* y u^* \in A^{\perp \perp}$, and so since $X^{\perp \perp}$
is an inner ideal we see that
$y = u u^* y u^* u \in X^{\perp \perp} \cap A = X$.
So $X = \, _pA_q = A_u$.  Since $(x d_m c_m d_n)$ converges
with $n$ for fixed $m$, for $x \in X$, we have $x d_m c_m u^* \in A$.
Taking the limit over $m$ shows that $x u^* \in A$.
By  Theorem \ref{mat1} we see that $u$ is a $*$-open tripotent.     \end{proof}

  \begin{xrem} Hereditary bimodules generalize the
hereditary subalgebras  of \cite{BHN} (called HSA's).  Indeed, HSA's are just the
case that $c_n = d_n$ in the characterization above.
If $c_n = d_n$ in the characterization above, then
$u = u^*$, so that $p = q$ in the proof above, and $X = \, _pA_p$ is
a HSA in $A$.  Conversely given a HSA with cai $(e_n)$,
the conditions above are met with $c_n = d_n = e_n$.
\end{xrem}

It is clear in the definition of a hereditary bimodule $X$  at the
start of this section that
$D \cong E$ completely isometrically as operator algebras,
and $D \cong X \cong Y$ canonically as operator algebras.
In fact these isomorphisms are easy to write down in terms
of the support tripotent $u$ for $X$:

\begin{proposition} \label{expis}  Let  $X$ be a hereditary bimodule
with a support tripotent $u$.
In  the notation at the start of this section we have
$D \cong E$ completely isometrically via the homomorphism $a \mapsto u^* a u$,  and
$X \cong Y$ completely isometrically via the linear isometry
$x \mapsto u^* x u^* \overset{{\rm def}}{=} x^\sharp$.   Similarly,
$D \cong X$ via right multiplication by $u$, and $D \cong Y$ via left
multiplication by $u^*$.
  \end{proposition}

\begin{proof}  That $D = \, _pA_p \cong \, _qA_q = E$ via this
map was done in the first paragraph of the proof of
Proposition \ref{mortri}; as were the
$D \cong X$ and $D \cong Y$ assertions.  Similar  arguments prove the
other assertion: If $X = A_u$ for a $*$-open tripotent $u \in \Delta(A^{**})$,
then the `dual module' $Y = \tilde{X}$ (see the definition
of a hereditary bimodule and the accompanying discussion
there), is simply $v X v$ where $v = u^*$.
Indeed as we saw in the proof of
Proposition \ref{mortri}, $X u^* = A_u u^* = A_{u u^*}$, and similarly
$u^* A_{u u^*} = A_{u^*} = Y$.
Thus we can define a complete isometry from $X$ onto $Y$ by $x \mapsto x^\sharp = v x v$.  The inverse
of this is of course the map $y \mapsto u y u$.
 \end{proof}

 \begin{xrem} (1) \  One should beware though in the
last result: the isomorphisms mentioned there are not unique, but depend on
the particular support tripotent $u$ chosen for $X$.  They are, however, determined for
example by the cone $\Rdb^+ {\mathfrak S}_{A_u}$ in Lemma \ref{ord} (iii).
(The nonuniqueness
of the support tripotent is discussed
a few paragraphs before Proposition  \ref{mortri}.)

(2) \  Consider the bijection between $*$-open tripotents $v$
and the elements $x$ and $y$ in the variants of item (ii) in Theorem \ref{varp},
described in  Proposition \ref{simt}.  One
may go further by saying that
 the $x$ and $y$ in Theorem \ref{varp} (ii) and its variants, are nothing but
the `square root of $a$ with respect to the tripotent' $v$.  That is, for example,
 via the canonical isomorphisms
$A_{a} \cong \overline{aAb}$ and $A_{a} \cong \overline{bAa}$ from Proposition
\ref{expis} (but with $u = v^*$), $x$ and $y$ both
correspond to $a^{\frac{1}{2}}$.  Of course $y=x^{\sharp}$ in the language of
Proposition \ref{expis}.  \end{xrem}

We can now give several characterizations of Pedersen equivalence in terms of hereditary bimodules:

\begin{corollary} \label{finma2}  If $a,b \in {\mathfrak S}_A$
then $a \sim_{A} b$ if and only if $\overline{aAb}$ is a
hereditary bimodule and there exists $x \in \overline{aAb}$
such that
$a = x \, x^\sharp, b=x^\sharp \, x$.   Here $x^\sharp = u^* x u^*$ is as in
Proposition {\rm \ref{expis}}, where
$u$ is a  $*$-open tripotent with $\overline{aAb} = A_u$.

If these hold then  $\overline{aAb}$ is indeed a  principal hereditary $A_a$-$A_b$-bimodule.
 \end{corollary}
\begin{proof}   ($\Rightarrow$) \ If $a \sim_{A} b$, and $u$ is as in
Proposition \ref{simt}, and  $v = u^*$,  then $A_u = \overline{aAb}$ by
Proposition \ref{aab}.  By Proposition \ref{mortri} this is a
hereditary bimodule.   The rest is clear from Remark (2) above.

($\Leftarrow$) \  Follows from Proposition \ref{simt},
since if $v = u^*$ then $$vav^{*} = v x v x v v^{*} = x^{\sharp} \, x = b .$$  
If these hold, then as we said  $\overline{aAb}$
is a hereditary bimodule.   We know from \cite{BRead} that
$\overline{A_a} = \, _pA_p$ and $\overline{A_b} = \, _qA_q$,
where $p = p_a, q = p_b$.
Also, $\overline{aAb} = \, _pA_q$ by Proposition \ref{aab},
and similarly, $\overline{bAa} = \, _qA_p$.
So $\overline{aAb}$ is a principal hereditary $A_a$-$A_b$-bimodule.
\end{proof}
\begin{corollary}   \label{finma}
If $a,b \in {\mathfrak S}_A$ and $\overline{aAb}$ is a
hereditary bimodule, with nets $c_n, d_n$ as in Theorem {\rm \ref{phii}},
then $a \sim_{A} b$ if and only if we can write $a = xy, b=yx$,  where
$y \in \overline{bAa}$ and $x \in \overline{aAb}$, and $yc_{n}-d_{n}x \to 0$ in norm.
\end{corollary}
\begin{proof}For the one direction, if $u$ is in the proof of
Theorem \ref{phii}, then $y u = u^* x$
because $yc_{n}-d_{n}x \to 0$.   So $y = x^{\sharp}$ and the rest
follows from Corollary \ref{finma2}.
For the other direction, if $a \sim_{A} b$ then
$a=xy,b=yx$ where $x = a^{1/2}v^* = v^* b^{1/2}$ and $y = v a^{1/2} =b^{1/2}
v$ for open tripotents $v, v^{*} \in A^{\perp\perp}$.
One can let $c_{n} = a^{1/2n} v^*, d_{n} = v a^{1/2n}$, and
the rest is easy.  \end{proof}

\begin{xrem}   (1) \ If $A = M_n$ then saying that
$\overline{aAb}$ is a hereditary bimodule is simply saying that
rank$(a) = {\rm rank}(b)$.  So we can view $\overline{aAb}$ being
a hereditary bimodule as a generalization of this rank condition.
 If we reduce the conditions in Corollary \ref{finma} to say: ($a = xy, b=yx$ where
 $\overline{aAb}$ is a hereditary bimodule,  $x \in \overline{aAb}$, and
$y \in \overline{bAa}$), then one obtains a
 relation that is in general strictly stronger
(by Corollary \ref{mst2}) than
the (variant of) `Blackadar equivalence' which we study later, and is
strictly weaker in general than Pedersen equivalence.

(2) \ Another characterization of Pedersen equivalence: If $a,b \in {\mathfrak S}_A$ then $a \sim_A b$ iff
  there exists a principal hereditary $A_a$-$A_b$-bimodule,
such that $a \oplus 0$ and $0 \oplus b$ in
the linking algebra of the bimodule correspond via the map $\rho$ mentioned in the first few paragraphs of
Section \ref{hbm},
to $c \oplus 0$ and $0 \oplus c$ for some element  $c \in C$.
For the one direction of this, in  Corollary \ref{finma2} we saw $\overline{aAb}$
is a principal hereditary bimodule.
If $v = u^*$ are as in the  proof of Corollary \ref{finma2}, with $v^* b v = a$,
and if $\rho^{-1}$ is the map in the first paragraph of the proof of
Proposition \ref{mortri},
then $\rho^{-1}(1 \oplus b) = 1 \oplus (v^* b v) = 1 \oplus a$.
For the other direction, by the proof of Proposition  \ref{mortri}, $\overline{aAb} = A_u$
for a $*$-open tripotent with $u u^* = p_a, u^* u = p_b$.  If $\rho$ is as in the second
paragraph of that proof,
and if $\rho(c \oplus 0) = a \oplus 0$, then
by hypothesis  $0  \oplus b = \rho(0 \oplus c)$.  However, $$\rho(0 \oplus c) =
\rho^{**}(E_{21} (c \oplus 0) E_{12}) =
\rho^{**}(E_{21}) (a \oplus 0) \rho^{**}(E_{12}) = 0 \oplus v a v^* .$$
Thus $b = v a v^*$. So (iv) in Theorem \ref{varp} holds.  \end{xrem}

\begin{definition}  \label{supptr}
If $u$ is a $*$-open tripotent in $A^{**}$ and $x \in {\mathfrak c}_{A_u}$,
where the latter denotes $\Rdb^+ {\mathfrak S}_{A_v}$,
then the $n$th roots of $x$ (with respect to the
Peirce product on $A_u$), converge weak* to an open projection
$s_u(x) \in B^{**}_2(u)$.   By Lemma \ref{ord},
$s_u(x)$ is a $*$-open tripotent in $A^{**}$, and $s_u(x) \leq u$ in the
ordering on tripotents.  We
 call $s_u(x)$ the $u$-{\em support tripotent} of $x$.  We also call $s_u(x)$
a {\em local support tripotent}.
\end{definition}

\begin{lemma} \label{hco}  Let $u$ be a $*$-open tripotent in $A^{**}$ and
$x \in {\mathfrak c}_{A_u}$.
The hereditary bimodule $A_s$ corresponding to the local support
tripotent $s= s_u(x)$ is
a principal hereditary bimodule, indeed may be written as
$\overline{aAb}$ for some $a, b \in {\mathfrak S}_{A}$ with $a \sim_A b$.
\end{lemma}

\begin{proof}  By \cite[Corollary 2.6]{BRead}, working within the algebra $A_u$,
and letting $s = s_u(x)$ and $v = u^*$, we have
$$\overline{x v A_u v x}  = \{ x \in A_u : x = s v x v s
= s s^* x s^* s \} = \{  x \in A : x =  s s^* x s^* s \} =
A_{s} .$$ Also $$x v = x v s v = x v s s^* s v =
x s^* s v = x s^* ,$$ and similarly $v x = s^* x$.
    Let $a = x s^*, b = s^* x = s^* a s$.   Since
$s$ is $*$-open in $A^{**}$, if $p = ss^*$
then right multiplication by $s^*$ is a completely
 isometric homomorphism from $A_s$ onto $_pA_p$,
which extends to a completely
 isometric homomorphism from $(A_s)^{\perp \perp}$
onto $p A^{**} p$.  It follows that $p_a = s s^*$.  Similarly $p_b =
s^* s$. Then $a$ and $b$ are in ${\mathfrak S}_{A}$, since for
example  $$\Vert 1 - a \Vert = \Vert u u^* - x u^* \Vert \leq \Vert
u - x  \Vert \leq 1.$$ Clearly $a s = x \in A$,  so $a \sim_A b$ by
Theorem \ref{varp} (iv)$^{\prime}$.

Clearly, $x v A_u v x = a A_u b \subset \overline{aAb}$.
To see the converse, first note that
$\overline{aA} = \overline{a^2 A}$.
Indeed, squaring the power series in Lemma \ref{pow} shows that $a \in \overline{a^{2}A}$ and the statement is thus clear.
Similarly $\overline{Ab} =\overline{Ab^2}$.   So
$$\overline{aAb} = \overline{aA} \overline{Ab} = \overline{a^2 A} \overline{Ab^2}
 \subset \overline{a A_u b},$$ the latter since $u u^*
x = x$ and $x u^* u = x$.  So $\overline{x u^* A_u u^* x}  = \overline{aAb}$, and
we are done.
 \end{proof}

 \begin{xrem} It is easy to see that $\overline{aAb}$ in the last proof,
actually equals $\overline{a A_{s_u(x)}}$, and
also equals $\overline{A_{s_u(x)}b}$.
\end{xrem}

\begin{theorem}  \label{anyh} \begin{enumerate}
 \item [{\rm (1)}]
Every hereditary bimodule in $A$ is
the closure of an increasing net of principal
hereditary bimodules of the form $\overline{aAb}$, where
$a, b \in {\mathfrak S}_A$ with $a \sim_A b$.
 \item [{\rm (2)}]   Any $*$-open tripotent $u$ in
$A^{**}$ is the limit of an increasing net of local support
tripotents.
\end{enumerate}
\end{theorem}

\begin{proof}  (1) \ Viewing $A_u$ as an operator algebra
with the Peirce product, by \cite[Theorem 2.15]{BRead}
$A_u$ is the closure of an increasing net $(D_t)$ of `peak principal' HSA's,
 each of the form $\overline{x_t A_u x_t}$ (with respect to the
Peirce product in $A_u$), for some $x_t \in {\mathfrak S}_{A_u}$.  As we saw in
Lemma \ref{hco}, we may write each $D_t$ as
$\overline{aAb}$ for some $a, b \in {\mathfrak S}_{A}$ with $a \sim_A b$.

 (2) \
 In the above, $s_u(x_t)$ forms an increasing net of open projections in
 $(A_u)^{\perp \perp}$ which converge weak* to $u$.  Hence by Lemma \ref{ord} they are local support tripotents,
and converging weak* to $u$.
\end{proof}

\begin{theorem}  \label{anys}  Any separable
hereditary bimodule in $A$ is a
principal
hereditary bimodule, that is of the form $\overline{aAb}$, where
$a, b \in {\mathfrak S}_A$ with $a \sim_A b$.  \end{theorem}

\begin{proof}    If $A_u$ is separable, then viewing it as
an operator algebra with respect to the
Peirce product,  we have by \cite[Corollary 2.18]{BRead} that
$A_u = \overline{x A_u x}$ with respect to the
Peirce product, for some $x \in {\mathfrak S}_{A_u}$.
By Lemma \ref{hco}, $A_u = \overline{aAb}$ for some $a, b \in {\mathfrak S}_{A}$ with $a \sim_A b$.
\end{proof}

 \begin{xrem} It follows as in \cite[Corollary 2.18]{BRead} that the
cai for a separable hereditary bimodule may be chosen to be countable, and
consist of mutually commuting elements in the Peirce product.  We are not
sure what implications this may have.  \end{xrem}

\begin{corollary} \label{vifa}  If $A$ is a separable operator algebra,
then every $*$-open tripotent in $A^{**}$ is a local support tripotent
$s_u(x)$ for some $x \in {\mathfrak S}_A$.
\end{corollary}  \begin{proof}  In the last proof,  we have $A_u = A_{s_u(x)}$,
with $s_u(x)$ a tripotent dominated by the $*$-open tripotent $u$.  By Lemma
\ref{ord}, $u = s_u(x)$.  \end{proof}

We say that a hereditary bimodule $D$ in $A$ is a
{\em hereditary subbimodule} of another hereditary bimodule $D'$ in $A$,
if $D \subset D'$, if the Peirce product of $D$ coincides with the restriction
of the Peirce product of $D'$, and $D D' D \subset D$ in the
Peirce product of $D'$.

\begin{proposition}  \label{ctsu}
If a hereditary bimodule $D$ in $A$ is generated
by a countable collection of principal hereditary subbimodules $D_n$
of $D$,  (so in the Peirce product on $D$ there is no
proper HSA containing all the $D_n$), then
$D$ is a principal hereditary  bimodule in $A$.
\end{proposition}

\begin{proof}  This follows from \cite[Theorem 2.16 (2)]{BRead} in the
same way that the last several results followed from matching results
in \cite{BRead}.
\end{proof}

For interests sake we give a generalization of
Theorem \ref{anys} to a slightly more general class of
inner ideals in $A$.  Similarly to our first
definition of a hereditary
bimodule at the start of this section, suppose that we have subalgebras $D, E$ of $A$,
and subspaces $X, Y$ of $A$, such that $(D,E,X,Y)$ is a {\em strong Morita equivalence context}
 in the sense of \cite{BMP}, and the
linking algebra of this context, namely
the subalgebra ${\mathcal L}$ of $M_2(A)$ with rows $C,X$ and $Y,D$, is
a HSA in $M_2(A)$.
Then we say that $X$ is
an {\em inner equivalence bimodule},
or an {\em inner equivalence} $D$-$E$-{\em bimodule}.
Note that again $D = XY, E = YX, X = DAE, Y = EAD$,
and $D$ and $E$ are HSA's and $X$ and $Y$ are inner ideals in $A$.  These follow by the same reasons as before.
A {\em principal inner equivalence bimodule},
or a {\em principal inner equivalence} $A_a$-$A_b$-{\em bimodule},
is the case when the context is
$(A_a,A_b,\overline{aAb}, \overline{bAa})$ for some $a, b \in {\mathfrak S}_A$.
Remark (1) before Theorem \ref{varp} essentially
says that if $a \sim_r b$
then $\overline{aAb}$ is a principal inner equivalence $A_a$-$A_b$-bimodule.

\begin{theorem}  \label{prie}  Any inner equivalence bimodule in $A$
which is also separable, is a principal inner equivalence bimodule.
\end{theorem}

\begin{proof} Let $(D,E,X,Y)$ is as above, and suppose that
the inner equivalence bimodule
$X$ is separable.   Then the ternary envelope $Z = {\mathcal T}(X)$ is
separable (see e.g.\  \cite{Bmuls,BLM}  for a
discussion of the ternary envelope, sometimes called the triple envelope).
By \cite[Lemma 10.2]{Bmuls},  the TRO $Z$ is
isomorphic to a canonical $C^*$-algebraic Morita
equivalence bimodule $Z$ containing $X$.   By facts in  e.g.\ p.\ 407 of  \cite{AGOH},
$Z^*$ contains $Y$, $Z Z^*$ contains $D$, and $Z^* Z$ contains $E$.
So $D, E, Y$ are each separable.  By \cite[Theorem 2.16 (1)]{BRead},
$D = A_a, E = A_b$ for some $a, b \in {\mathfrak S}_A$.
Thus $X = D A E \subset \overline{aAb}$, and conversely
$\overline{aAb} \subset  D A E =  X$ since $a \in D, b \in E$.  So $X = \overline{aAb}$,
and similarly $Y = \overline{bAa}$.
        \end{proof}

\section{Blackadar equivalence}

For $a,b \in  {\mathfrak S}_A$ we define
$a \cong b$ if $s(a) = s(b)$.  By the theory in
\cite[Section 2]{BHN} this is  equivalent to $\overline{aA} = \overline{bA}$, and also
 equivalent to $\overline{aAa} = \overline{bAb}$.

If $a, b \in {\mathfrak c}_A$ we define $a \sim_s b$ if there exist $a',b' \in
 {\mathfrak S}_A$,
with $a \cong a', a' \sim_A b'$, and $b' \cong b$.
It will follow from the next result that this is an equivalence relation.
We say that $a$ and $b$ are {\em Blackadar equivalent in}
$A$ if $a \sim_s b$.

\begin{theorem}  \label{isaq}  {\rm (Cf.\ \cite{Lin}, \cite[Proposition 4.3]{ORT}.)}  For $a,b \in   {\mathfrak c}_A$ TFAE:
\begin{enumerate} \item  [(i)] $a \sim_s b$.
\item   [(ii)] $\overline{aA} \cong \overline{bA}$
completely isometrically via a right $A$-module map.
\item  [(iii)] $p_a \sim_{PZ,A} p_b$.
\item   [(iv)] There exists $b' \in  {\mathfrak c}_A$,
with $a \sim_A b'$ and $b' \cong b$.
\item   [(v)] There exists $a' \in
  {\mathfrak c}_A$, with $a \cong a' \sim_A b$.
\item   [(vi)]  $\overline{aAb}$
is  a
principal hereditary  $A_a$-$A_b$-bimodule.  \end{enumerate} \end{theorem}

\begin{proof}  We may assume that $a,b \in   {\mathfrak S}_A$.

(iv) $\Rightarrow$ (vi) \ There exists $f \in  {\mathfrak c}_A$
with $a \sim_A f$ and $f \cong b$.   By Proposition \ref{finma2},  $\overline{aAf}$
is a principal hereditary bimodule,
but
 $\overline{fAf} = \overline{bAb}$.

(vi) $\Rightarrow$ (iii) \ By Proposition  \ref{mortri} we have
$X = \overline{aAb} = A_u$ for a $*$-open tripotent $u \in A^{\perp\perp}$.
Note that $p_a x = x$ for all $x \in X$, so that $p_a u = u$.
Hence $p_{a}u u^{*} = u u^{*}$, so $u u^{*} \leq p_a$.   Similarly,
$u^{*} u \leq p_{b}$.
However since $a \in A_a = X Y$, we have $u u^{*} a = a$,
so that $u u^{*} \leq p_a$.  Hence $u u^{*} = p_a$ and similarly
$u^{*} u = p_b$.   Of course $u a \in A$ by definition of a $*$-open tripotent.

(iv) $\Rightarrow$ (i) \ Trivial.

(ii) $\Rightarrow$ (iv) \   If $\Phi : \overline{aA} \to \overline{bA}$
is a surjective completely isometric right $A$-module map,
then by \cite[Corollary 3.7]{AGOH} there exists a
surjective completely isometric left $A$-module map  $\Psi
: \overline{Aa} \cong \overline{Ab}$  such that $\Psi(y) \Phi(x) = yx$
for all $x \in aA, y \in Aa$.

We obtain using \cite[Theorem 6.8]{AGOH}, a $C^*$-module isomorphism
$$B \otimes_{hA} \overline{Aa} \to B \otimes_{hA} \overline{Ab}.$$
Now the multiplication map $m: B \otimes_{hA} \overline{Aa} \to
\overline{Ba}$ has an asymptotic contractive left inverse
 $\theta_n : z \mapsto z \otimes a^{\frac{1}{n}}$.
Indeed $$\theta_n(m(b \otimes x)) = bx \otimes a^{\frac{1}{n}} = b
\otimes x a^{\frac{1}{n}} \to b \otimes x.$$ for all $b \in B, x \in
\overline{Aa}$.  So $m$ is isometric, and so $\Psi$ extends to a
$C^*$-module map $\overline{Ba} \to  \overline{Bb}$ which is unitary
onto its range.  It follows from e.g.\ \cite[Corollary 8.1.8]{BLM}
that $y^* y = (a^{\frac{1}{2}})^* a^{\frac{1}{2}}$, where $y =
\Phi(a^{\frac{1}{2}})$. Similarly, $x x^* = a^{\frac{1}{2}}
(a^{\frac{1}{2}})^*$ if $x = \Psi(a^{\frac{1}{2}})$, and $x y =
\Psi(a^{\frac{1}{2}}) \Phi(a^{\frac{1}{2}}) = a^{\frac{1}{2}}
a^{\frac{1}{2}}) = a$.  By Theorem \ref{vpr}, $b' = yx \in
{\mathfrak c}_A$ and $a \sim_A b'$.  Finally, $\overline{b' A} =
\overline{yxA} \subset \overline{yA}$. On the other hand, $y =
\lim_n \, y a^{\frac{1}{n}} \in y \overline{aA} \subset
\overline{yxA}$, since $a = xy$.  So $$\overline{b' A} = \overline{y
A} = \overline{\Phi(a^{\frac{1}{2}}) A} =
 \Phi(\overline{a^{\frac{1}{2}} A}) =
\overline{bA}.$$   Thus $a \sim_A b' \cong b$.

(i) $\Rightarrow$  (iii) \ If $a \sim_A b$
and $v$ is as in the proof of Theorem \ref{varp}
above, then $v, v^* \in A^{\perp \perp}$ as we said
in Theorem \ref{varp}.   Thus $p_a \sim_{PZ,A} p_b$ via
$v$ (for example, $v^* b = v^* d d = xd \in A$
so $v^* B_q = v^* \overline{bBb} \subset B$).  Since $p_a = p_{a'}$ if $a \cong a'$
(iii) is now clear.

(iii) $\Rightarrow$  (ii) \   If  $p_a \sim_{\rm PZ,A} p_b$ via
a partial isometry $v$, let $\Phi(x) = vx$ and $\Psi(x) =
x v^*$.  These are the desired
completely isometric module maps
in (ii) (note that $\Phi(aA) = vaA \subset v A_p A \in A \cap
q B^{**} = \overline{bA}$,
 and similarly $\Psi(bA) \subset aA$).

The proof of the equivalence with (v) is similar to
the proof of the equivalence with (iv). \end{proof}

\begin{xrem}  For algebras without any nontrivial r-ideals all elements
of ${\mathfrak c}_A$ are Blackadar equivalent to each other of course.  Such algebras are
discussed in \cite{BRead}.  \end{xrem}

\begin{corollary} \label{mst2}  If $a, b \in  {\mathfrak c}_A$,
and $X = \overline{aAb}$ and $Y = \overline{bAa}$ then TFAE:
\begin{enumerate} \item  [(i)]  $p_a \sim_{PZ,A} p_b$.
 \item  [(ii)]  There
are nets of contractions $c_n \in X$, and $d_n \in A$ satisfying
{\rm (2)} of Theorem {\rm \ref{phii}}, and also
$c_n d_n a \to a$ and $b d_n c_n \to b$ in norm.
\item   [(iii)] $X$ is a hereditary  bimodule
and also  $a \in XY, b \in YX$.
  \end{enumerate} \end{corollary}

\begin{proof}  (ii) $\Rightarrow$ (i) \ These conditions guarantee by Theorem \ref{phii}
that $X$ is a hereditary bimodule, and
$X = A_u$ for a $*$-open tripotent $u \in A^{\perp\perp}$.
As in the proof of Theorem \ref{isaq},
$u u^{*} \leq p_{a}$ and $u^{*} u \leq p_b$.
However, the last condition in (ii), together with
a part of the proof of Theorem \ref{phii}, guarantees that $u u^{*}a =a$ and $bu^{*} u =b$.  So
$u u^{*} = p_{a}$ and $u^{*} u = p_{b}$.

(i) $\Rightarrow$ (iii) \   Follows from Theorem \ref{isaq}.

(iii) $\Rightarrow$ (ii) \ Follows from Theorem \ref{phii}.
 \end{proof}

We now give some other conditions that imply Blackadar equivalence in a $C^*$-algebra.

\begin{theorem} \label{m1}
 Suppose that $a, b \in {\mathfrak c}_A$ for an operator algebra $A$. If
 $a = xy$ and $b = yx$ for some $x,y \in A$, then $p_a \sim_{PZ} p_b$
in $B^{**}$.
\end{theorem}

\begin{proof}  Without loss of generality, by dividing by a suitable scalar,
we may assume that $a, b \in {\mathfrak S}_A$.
Under these conditions $a x = xy x = xb$, and
hence $a^{\frac{1}{n}}x = x b^{\frac{1}{n}}$ for all $n \in \Ndb$ (as in Remark (2)
before Theorem \ref{varp}).  In the limit, $px = xq$ and $x^{*}p = qx^{*}$,
where $p = p_a$ and $q = p_b$.  It follows that $p x x^* x = x x^* x q$, and similarly
for fivefold and sevenfold products of this type, etc.
Thus $p r(x) = r(x)q$, where $r(x)$ is the range tripotent of $x$ in $Z^{**}$,
because $r(x)$ is a weak* limit of linear combinations
of such `odd products' (see e.g.\  the proof of
\cite[Lemma 3.3]{BN0}). Since $p$ is a limit of functions of $xy$, clearly $r(x)r(x)^{*}p = p$. Thus $w=p \, r(x)$
is a tripotent and it is easy to see that $ww^{*}=p$ and $w^{*}w=q$.
Note that $w^* a = r(x)^* xy = |x|y \in B$.  So $w^* \, _pB_p =
w^* \overline{aBa} \subset B$.
By \cite[Lemma 1.3]{PZ} we have
$p_a \sim_{PZ} p_b$.
 \end{proof}

It is not necessarily
true that if  $a = xy$ and $b = yx$ for some $x,y \in A$, then $a$ and $b$ are Blackadar equivalent in $A$,
if $A$ is nonselfadjoint.
For a counterexample, let $R$ be an invertible operator on a Hilbert space $H$, and let $A$ be the span in
$M_2(B(H))$  of
$a = I_H \oplus 0, b = 0 \oplus I_H, x = E_{12} \otimes R, y = E_{21} \otimes R^{-1}$.  Then  $a = xy$ and $b = yx$,
but $a$ and $b$ are not Blackadar equivalent in $A$ if $R$ is chosen appropriately.

This counterexample also shows that the conditions considered in (2) of the next result also do not characterize
 Blackadar equivalence in $A$ if $A$ is nonselfadjoint.

\begin{proposition}  Let $A$ be an approximately unital operator algebra, suppose that $a, b  \in  {\mathfrak c}_A$, and let $B$ be a $C^*$-algebra containing $A$ as
usual.  Then \begin{enumerate}
 \item Suppose that  for some $x \in A$ we have $\overline{aA} =  \overline{x  A}$, and $\overline{Ab} =  \overline{Ax}$.  Then $a$ and $b$ are Blackadar equivalent in $B$.
 \item  If $a$ and $b$ are Blackadar equivalent in $A$
then there exists an $x \in A$ such that
 $\overline{aA} =  \overline{x  A}$, and $\overline{Ab} =  \overline{Ax}$.
\item $a$ and $b$ are Blackadar equivalent in $B$ if and only if there exists an $x \in B$ such that
 $\overline{aA} =  \overline{x  A}$, and $\overline{Ab} =  \overline{Ax}$.  \end{enumerate} \end{proposition}

\begin{proof}  (1) \ Let $p = p_a, q = p_b$.  If these conditions hold then $x \in \overline{aA} \cap
\overline{Ab}$, so
 $px = pxq
= xq$.  Since $a  \in
\overline{xA}$ we have $r(x)r(x)^{*}p = p$, and the proof of (1) then follows as in the previous theorem.

(2) \   These follow  for example  by the same reasoning used in the proof of Theorem \ref{isaq} to show that
$\overline{yA} =  \overline{b'A} =\overline{bA}$.

(3) \  This is now obvious.  \end{proof}

 \begin{xrem} (1) \  If $A$ is a $C^*$-algebra then the last result gives an alternative characterization of
Peligrad-Zsid\'o/Blackadar equivalence.
One may ask what is the least restrictive extra condition we can add to these to characterize our variant of
Blackadar equivalence in a
more general operator algebra.   For example, one such extra condition is  that $\overline{aAb}$ is
a hereditary bimodule (see Corollary \ref{mst2} (iii)).

(2) \ Another relation equivalent to the one considered in the last result, is that
  $\overline{A_a \, x}= \overline{x \, A_b}$ for some $x \in A$, and $a \in
\overline{xA}$ and $b \in \overline{Bx}$.  This also characterizes
Peligrad-Zsid\'o/Blackadar equivalence in a $C^*$-algebra.
\end{xrem}

In the remainder of this section, we turn to the topic of subequivalence.

For any $a,b  \in {\mathfrak c}_A$ we have as in
4.4 in \cite{ORT},  that
$a \in A_b$ iff $A_a \subset A_b$ iff $\overline{aA} \subset \overline{bA}$ iff
$p_a \leq p_b$.  Also ($a \in A_b$ and $b \in A_a$) iff
 $a \cong b$ iff $\overline{aA} = \overline{bA}$ iff
$p_a = p_b$.

\begin{lemma} \label{haver}  Suppose that $a \in {\mathfrak c}_A$, and
that $\Phi, \Psi$ are a pair of completely contractive
$A$-module maps from $\overline{aA}$ and $\overline{Aa}$ into
$A$, such that $\Psi(x) \Phi(y) = xy$ for all $x \in \overline{Aa},
y \in \overline{aA}$.   Then $\Phi, \Psi$ are completely isometric, and
there exists $b \in  {\mathfrak c}_A$ with $a \sim_A b$
and  ${\rm Ran}(\Phi) = \overline{bA}$ and ${\rm Ran}(\Psi)
= \overline{Ab}$.   \end{lemma} \begin{proof}
Let $F, E$ be the ranges of these two maps,
which are respectively a right and a left ideal in $A$.
By applying the $\otimes_{hA}$ tensor product with
$B$ as in the proof that (ii) $\Rightarrow$ (iv)
in Theorem \ref{isaq}  above, we may extend  $\Phi, \Psi$ to
contractive $B$-module maps from $\overline{aB} $ and $\overline{Ba}$ into
$B$.  If $\Phi(c) = y, \Psi(c) = x$, where $c = a^{\frac{1}{2}}$,
then $y^* y \leq c^* c$  by e.g.\ 8.1.5 in \cite{BLM}, and
similarly $x x^* \leq c c^*$.    By basic operator theory, $y = Sc, x = cR$
for contractions $S, R$.  By the proof that (ii)$^{\prime \prime \prime}$
implies  (ii)$^{\prime}$ in Theorem \ref{varp}, we actually have
 $|y| = |c|, |x^*| = |c^*|$.   Thus $\Phi(cz)^* \Phi(cw) = z^* y^* y w = (cz)^* (cw)$,
for all $z, w \in A$, which forces $\Phi$ to be a complete isometry.  Similarly
$\Psi$ is  a complete isometry.
 Setting $b = yx$, then $b \in
 {\mathfrak c}_A$ by Theorem \ref{vpr}, and
$a \sim_A b$.
Also, $bA = yx A \subset \overline{yA}$, and conversely,
$y a^{\frac{1}{n}} \in y  \overline{aA} \subset
 \overline{bxA} \subset \overline{bA}$,
so $y \in \overline{bA}$ and $\overline{yA} \subset \overline{bA}$.
Thus $\overline{yA} = \overline{bA}$.   Hence
$F = \Phi(\overline{aA}) = \Phi(\overline{cA}) = \overline{yA} = \overline{bA}$.
Similarly, $E = \overline{Ab}$.  \end{proof}

 If $a, b \in  {\mathfrak c}_A$,
we define $a \precsim_{s} b$ if there exists $b' \in A_b$ such that
$a \sim_s b'$.  This is clearly equivalent to:
there exists $b' \in A_b$ such that
$a \sim_A b'$.   We will call this {\em Blackadar comparison in}
$A$.   If $p, q$ are open projections
in $A^{**}$ we say that $p \precsim_{PZ,A} q$ if there is
an open projection $q' \leq q$ in  $A^{**}$ with
$p \sim_{PZ,A} q'$.  We will call this  {\em Peligrad-Zsid\'o subequivalence in}
$A^{**}$.

The following is the version of \cite[Proposition 4.6]{ORT} in our setting
(see also \cite{Lin}):

\begin{proposition} \label{46}  If  $a,b  \in {\mathfrak c}_A$, TFAE:
\begin{enumerate} \item   [(i)] $a \precsim_{s} b$.
\item  [(ii)]  $p_a \precsim_{PZ,A} p_b$.
\item  [(iii)]  There exist
a pair of completely contractive
$A$-module maps $\Phi : \overline{aA} \to \overline{bA}$ and $\Psi
: \overline{Aa} \to \overline{Ab}$, such that $\Psi(x) \Phi(y) = xy$ for all $x \in \overline{Aa},
y \in \overline{aA}$.
\end{enumerate}
\end{proposition}

\begin{proof}   (iii) $\Rightarrow$ (i) \ By Lemma \ref{haver}, (iii) implies
that there exist $x, y \in {\rm Ball}(A)$ with $a = xy, b' = yx \in  {\mathfrak c}_A$.
In our case, $b' = yx \in A_b$ too.

(i) $\Rightarrow$ (ii) \  By Theorem \ref{isaq}, $p_a \sim_{PZ,A} p_{b'}$, and
clearly $p_{b'} \leq p_b$.

(ii) $\Rightarrow$ (iii) \  If we have an open projection $q' \leq p_b$ in  $A^{**}$ with
$p_a \sim_{PZ, A} q'$ via a tripotent $v$, then by Proposition \ref{isu2} we have $q' =
v v^* = p_{b'}$ where
$b' = v a v^*$.  Left multiplication by
$v$ is a completely isometric module map from $\overline{aA}$ into
$A$.  Note that $p_b va = p_b v v^* va = p_b p_{b'} va
= p_{b'} va = va$.   So the map maps into $_{p_b}A = \overline{bA}$.
Similarly, right multiplication by
$v^*$ is a completely isometric module map from $\overline{Aa}$ into
$\overline{Ab}$.   The rest is clear.   \end{proof}

As in \cite[Proposition 4.6]{ORT}, these give $x, y \in {\rm Ball}(A)$
such that $\overline{aA} = \overline{xy A}$ and $\overline{yxA} \subset \overline{bA}$,
and so on.  Here $yx \in  {\mathfrak c}_A$.

 \begin{xrem} In contrast to Theorem  \ref{isaq} (ii), one needs the statement
about $\Psi$ in (iii) of Proposition \ref{46}.   It is not automatic,
as one may see by considering for example multiplication by $z$ on the disk algebra.
  \end{xrem}

\section{Appendix: Some related results for TROs}

We include the following results about TROs for the reason that,
 as we claimed in Remark (2) at the start of
Section \ref{optoo}, many of our results generalize easily to the
case that the operator algebra $A$ is replaced by a strong Morita equivalence bimodule
in the sense of \cite{BMP}.   The results below show how this can proceed in the case of TROs.

If $v$ is a tripotent in the second dual
of a TRO
$Z$ such that $p = v^* v$ and $q = v v^*$ are open projections with respect to
$Z^* Z$ and $Z Z^*$ respectively, and if $v (Z Z^*)_p \subset Z, v^* (Z Z^*)_q \subset Z^*$,
then we say that $v$ implements a {\em TRO Peligrad-Zsid\'o equivalence} of open projections.
This is of course equivalent to $v$, thought of as in the second dual
of the linking $C^*$-algebra of $Z$, implementing a Peligrad-Zsid\'o equivalence
in the original sense between $0 \oplus p$ and $q \oplus 0$.  So again we need not say
$v^* (Z Z^*)_q \subset Z^*$ in this definition.  The following proof shows that
the two `inclusion conditions' are  equivalent to:
$v^* \, (_qZ_p)  \subset Z^* Z$ and $v \, (_pZ^*_q) \subset Z Z^*$.

\medskip

The following (which we should have seen in \cite{BN0}),
is the appropriate variant of the equivalence (i) $\Leftrightarrow$ (vi) in \cite{BW}:

\begin{proposition} \label{pzis}  If $v$ is a tripotent in the second dual
of a TRO
$Z$, then $v$ implements a TRO Peligrad-Zsid\'o equivalence of open projections
iff $v$ is an open tripotent in the sense of {\rm \cite{BN0}}.
\end{proposition}

\begin{proof}   ($\Rightarrow$) \
Under these conditions, note that $$D = Z(v) = Z^{**}_2(v) \cap Z =
\{ b \in Z : q b p = b \}$$ is an inner ideal in $Z$.  Here $p = v^*
v, q = v v^*$. It is also a $C^*$-subalgebra of $Z^{**}(v)$. Indeed
$$v^* D = v^* D D^* D \subset Z^* D \subset Z^* Z$$ and $$v D^* = v D^*
D D^* \subset Z D^* \subset Z Z^*,$$ hence $D v^* D \subset Z \cap
Z^{**}(v) = D$, and $v D^* v = vD^* D D^* v \subset Z Z^* D Z^* Z
\subset Z$, so $D v^* D \subset Z \cap Z^{**}(v) = D$.   If $B = Z^*
Z$ then $B_p = \{ b \in B : p  b p = b \}$ is a HSA with support
projection $p$, and $v^* D = B_p$ (clearly $v^* D \subset B_p$ and
conversely if $b \in B_p$ then $v b \subset Z \cap Z^{**}_2(v) = D$,
so $b \in v^* D$).  So left multiplication $L_{v^*}$ by $v^*$ is a
linear complete isometry from $D$ onto $B_p$, and
 is a ternary isomorphism, and indeed
it is a $*$-isomorphism with respect to the $v$-product on $D$.  It follows that
$D^* D = B_p$, and  similarly $D D^* = (Z Z^*)_q$, although these are clear
directly:    $D^* D = D^* v v^* D = B_p^* B_p = B_p$.

Suppose that $x_t$ is a positive cai in the $C^*$-algebra $D$ with the $v$-product,
with weak* limit $w$.
This will be an open tripotent in the sense of {\rm \cite{BN0}}
dominated by $v$.  Then $v^* x_t$ is a positive cai in
$B_p$, so that in the limit we have $v^* w = p$.   Since $w p = w$ and $v p = v$,
it follows from a basic operator theory fact that $v = w$.  So $v$ is open
in the sense of {\rm \cite{BN0}}.

($\Leftarrow$) \ If $v$ is open, then by the remark after Corollary 2.11
in \cite{BN0}, $p$ and $q$ are open.  If $x \in Z(v)$ then
$v^* v x^* x = x^* x \subset Z^* Z$.  Since $v x^* x$ is the generic positive
element in the $C^*$-algebra $Z(v)$, it follows that
$v^* Z(v) \subset Z^* Z$.  Similarly, $Z(v) v^*  \subset Z Z^*$.
We note that the HSA $Z(v) Z(v)^*$ of $Z Z^*$ has support projection $q$,
by the first part of the proof, and similarly
 $Z(v)^* Z(v)$ has support projection $p$.
So $v^* \,  (Z Z^*)_q = v^* Z(v) Z(v)^*  \subset Z^* Z \, Z(v)^*\subset Z^*$,
and similarly  $v \, (Z Z^*)_p \subset Z$.    \end{proof}

\begin{proposition} \label{isu}  Let $Z$ be a TRO,
set $A = Z Z^*, B = Z^* Z$,
and let $v$ be an open tripotent in $Z^{**}$,
and suppose that $p = v^* v$ is the support projection for some $b \in {\mathfrak S}_B$.
Then $q = v v^*$ is the support projection for some $a \in {\mathfrak S}_A$,
and in this case there exist $x, y \in {\rm Ball}(Z)$ with
$x \in \, _qZ_p, y \in \, _pZ_q$, and $x y = a, yx = b$.
Moreover, $\overline{aAa}$ and $\overline{bBb}$
are $*$-isomorphic $C^*$-algebras. They are also
 strongly Morita equivalent  via the equivalence bimodule
$\overline{a Z b} = \overline{bZ^*a}^*$, and this bimodule
is ternary isomorphic to $\overline{bBb}$.
  \end{proposition} \begin{proof}    As in Proposition 3.3 of \cite{ORT}, let $a = v b v^*$.  Note that
$a = v b^{\frac{1}{2}} b^{\frac{1}{2}}
v^* \in v B_p (v B_p)^* \subset A_q$.  Also,
$a \in {\mathfrak S}_A$ by Lemma \ref{cofr},
and
 $(v b v^*)^{\frac{1}{n}} =  v b^{\frac{1}{n}} v^*$
 by Lemma \ref{vpow}.  Then
$$p_a = \lim_n \, (v b v^*)^{\frac{1}{n}} =
\lim_n \, v b^{\frac{1}{n}} v^* = v p_b v^* = q ,$$
where these are weak* limits.   Hence $q$ is the support projection of $a$.  If $x = v b^{\frac{1}{2}}$ and $y =
b^{\frac{1}{2}} v^*$, then $x y = a, y x = b$ since $b^{\frac{1}{2}} \in B_p$.

Continuing as in \cite{BRead},  $\overline{bBb} =  B_p$ since $p =
s(b)$, and similarly $\overline{aAa} =  A_q$.  These are
$C^*$-algebras in this case. Also $a Z b \subset \, _qZ_p$, and,
conversely, we have $_qZ_p = A_q \, _qZ_p B_p \subset
\overline{aZb}$. So $\overline{a Z b}
 = \, _qZ_p$, and similarly $\overline{bZ^*a} =  \, _pZ^*_q$.
So we have a $C^*$-algebraic Morita equivalence between $\overline{aAa}$
and $\overline{bBb}$, implemented by the equivalence bimodule
$\overline{a Z b} = \overline{bZ^*a}^*$.

The last assertions follow from the first paragraph of Proposition
\ref{pzis}. Note that $z \mapsto v z v^*$ is a $*$-isomorphism
from $A_p$ onto $B_q$, with inverse the map $z \mapsto v^* z v$.      \end{proof}

\begin{proposition} \label{sep}
Any separable inner ideal $D$ in a TRO
$Z$ is of the form
 $\overline{a Z b}$ for $a \in (Z Z^*)_+, b \in (Z^* Z)_+$.  If $D$ is also ternary isomorphic
to a $C^*$-algebra then this can be done with $p_a \sim_{PZ} p_b$.  \end{proposition}

\begin{proof}  Let $A = Z Z^*, B = Z^* Z$.
If $a \in A_+, b \in B_+$ then
$\overline{a Z b}$ is an inner ideal.   Conversely, if $D$ is
an inner ideal in $Z$, then $D D^*$ is a HSA in $A$, which is separable
if $D$ is separable.  Hence $D D^* = \overline{a A a}$ for some $a \in
A_+$.  Similarly, $D^* D  = \overline{bBb}$  for some $b \in B_+$.
Thus $D = D D^* D D^* D = \overline{a Z b}$ (since
$a Z b \subset D (D^* Z D^*) D \subset D$ because
$D$ is an inner ideal).
If $D$ is also ternary isomorphic
to a $C^*$-algebra $C$ then the identity in $C^{**}$ corresponds
to an open tripotent $u \in D^{**}$ by \cite[Proposition 3.5]{BN0},
and it is easy to argue that $u u^*$ (resp.\ $u^* u$) is
the support projection of
$D D^*$ (resp.\ $D^* D$).  So $u u^*$ (resp.\ $u^* u$) is
the support projection of $a$ (resp.\ $b$).
\end{proof}

Acknowledgements:  
The second author was supported by Denison University.

\end{document}